\documentclass{compositio}
\usepackage{tikz-cd, amsmath, amsthm, amssymb, amsfonts}

\newtheorem{theorem}{Theorem}[section]
\newtheorem{lemma}[theorem]{Lemma}
\newtheorem{corollary}[theorem]{Corollary}
\newtheorem{proposition}[theorem]{Proposition}

\theoremstyle{definition}
\newtheorem{definition}[theorem]{Definition}

\newtheorem{remark}[theorem]{Remark}

\newcommand{\dval}{\mathrm{dval}}
\newcommand{\val}{\mathrm{val}}
\newcommand{\SRL}{\mathrm{RecMaps}}

\newcommand{\ord}{\mathrm{ord}}
\newcommand{\spec}{\mathrm{Spec}}
\newcommand{\Hom}{\mathrm{Hom}}
\newcommand{\ts}[1][3]{\partial^{(#1)}}

\def\t{\times}
\def\mc{\mathcal}
\def\mb{\mathbb}
\def\wt{\widetilde}
\def\wh{\widehat}
\def\bs{\backslash}
\def\CD{\mathcal P_2}
\def \C{\mathbb C}
\def\Q{\mathbb Q}
\def \N{\mathcal N}
\def \R{\mathbb R}
\def \Z{\mathbb Z}

\def\emb{\hookrightarrow}
\def \oF{\overline F}
\def \P{\mathcal B_2}
\def \ext{\text{ext}}
\def \nn{{\nu'|\nu}}

\def \oL{\overline L}
\def \oM{\overline M}
\def \oK{\overline K}
\def \ol{\overline}
\def\Set{\textbf{Set}}
\def\F{\mathbf{Fields}}

\begin{document}

\title{Chow dilogarithm and strong Suslin reciprocity law}
\author{Vasily Bolbachan}
\email{vbolbachan@gmail.com}
\address{Faculty of Mathematics, National Research University Higher School of Ecnomics, Russian Federation, Usacheva str., 6, Moscow, 119048; HSE-Skoltech International Laboratory of Representation
Theory and Mathematical Physics, Usacheva str., 6, Moscow, 119048}
%

%
%\dedication{to Alexander Goncharov.}


\thanks{This paper was partially supported by the Basic Research Program at the HSE University and by the Moebius Contest Foundation for Young Scientists}

\begin{abstract}
We prove a conjecture of A. Goncharov concerning strong Suslin reciprocity
law. The main idea of the proof is the construction of the norm map on so-called
lifted reciprocity maps. This construction is similar to the construction
of the norm map on Milnor $K$-theory. As an application, we express  Chow dilogarithm in terms of Bloch-Wigner dilogarithm. Also, we obtain a new reciprocity law for four rational functions on an arbitrary algebraic surface with values in the pre-Bloch group.
\end{abstract}

\maketitle
%\tableofcontents

\vspace*{6pt}\tableofcontents  % for this guide only.

\section{Introduction}
\label{sec:intro}
Everywhere we work over $\Q$. So any abelian group is supposed to be tensored by $\Q$. For example, when we write $\Lambda^2 k^\t$ this actually means $\left(\Lambda^2 k^\t\right)\otimes \Q$. All exterior powers and tensor products are over $\Q$. 

Let $k$ be a field and $X$ be a smooth projective curve over $k$. For a field $F$ denote by $K_n(F)$ the $n$-th algebraic $K$-theory of $F$. 
For any closed point $z\in X$ one can define the residue map $\widehat\partial_z\colon K_n(k(X))\to K_{n-1}(k(z))$ (we use this notation to distinguish this map from the residue map on polylogarithmic complexes which will be defined below), where $k(z)$ is the residue field of the point $z$ (see \cite[V.5]{KBook}). Denote by $tr_{k(z)/k}$ the push-forward map $K_{n-1}(k(z))\to K_{n-1}(k)$ associated to the natural projection $\spec( k(z))\to \spec (k)$. It follows from the basic properties of algebraic $K$-theory that for any $a\in K_n(k(X))$ and all but finitely many $z\in X$, we have $\widehat\partial_z(a)=0$ and moreover the following sum is equal to zero:
$$\sum\limits_{z\in X^{(1)}}tr_{k(z)/k}\circ \widehat\partial_z(a)=0.$$
In this formula $X^{(1)}$ denotes the set of closed points of the curve $X$.
 
On the other hand, for any field $F$, A. Goncharov \cite{goncharov1995geometry} defined so-called polylogarithmic complexes $\Gamma(F,n), n\in \mathbb N$ and  conjectured that these complexes compute the  graded pieces of the algebraic $K$-theory of $F$. More precisely the cohomology $H^i(\Gamma(F, n)\otimes \mathbb Q)$ should be isomorphic to $gr_\gamma^n K_{2n-i}(F)\otimes \mathbb Q$. Here $gr_\gamma^nK_{2n-i}(F)\otimes \mathbb Q$ is the associated graded space with respect to $\gamma$-filtration (see for example \cite{KBook}). 

The complex $\Gamma(F, n)$ looks as follows:
$$\Gamma(F,n)\colon \mathcal B_n(F)\xrightarrow{\delta_n} \mathcal B_{n-1}(F)\otimes F^\times\xrightarrow{\delta_n}\dots\xrightarrow{\delta_n}\mathcal B_2(F)\otimes \Lambda^{n-2}F^\times\xrightarrow{\delta_n}\Lambda^n F^\times.$$

This complex is concentrated in degrees $[1,n]$. The group $\mathcal B_n(F)$ is the quotient of the free abelian group generated by symbols $\{x\}_n, x\in \mathbb P^1(F)$ by some explicitly defined subgroup $\mathcal R_n(F)$ (see \cite{gon94}).  In the next section we will present the generators for the group $\mathcal R_2(F)$. The differential is defined as follows: $\delta_n(\{x\}_k\otimes y_{k+1}\wedge \dots \wedge y_n)=\{x\}_{k-1}\otimes x\wedge y_{k+1}\wedge \dots \wedge y_n$ for $k>2$ and $\delta_n(\{x\}_2\otimes y_3\wedge\dots y_n)=x\wedge (1-x)\wedge y_3\wedge \dots \wedge y_n$.

Let us assume that the field $k$ is algebraically closed. In this case  A. Goncharov constructed the morphism of complexes  $\partial_z^{(n)}\colon \Gamma(F,n)\to \Gamma(k,n-1)[-1]$ which should correspond to the residue map on the algebraic $K$-theory. So, it is natural to suppose that there is a homotopy between the map $\sum\limits_{z\in X^{(1)}}\partial_z^{(n)}$ and the zero map. In this paper we will deal only with the case $n=3$. In this case the existence of such a homotopy was proved in \cite{rudenko_2021}.  

It turns out that this story is connected with so-called Chow dilogarithm defined by A. Goncharov in \cite{goncharov2005polylogarithms}. For any smooth projective curve $X$ over $\mathbb C$ and three non-zero rational functions $f_1, f_2, f_3$ on $X$, the Chow dilogarithm $\mc P_2(X; f_1, f_2, f_3)$ is defined by the formula 
$$\mc P_2(X; f_1, f_2, f_3)=(2\pi i)^{-1}\int_{X(\C)}r_2(f_1, f_2, f_3),$$
where
$$r_2(f_1, f_2, f_3)=\dfrac 16\sum\limits_{\sigma\in S_3}sgn(\sigma)\wt r_2(f_{\sigma_1}, f_{\sigma_2}, f_{\sigma_3})$$

$$\wt r_2(f_1, f_2, f_3)=\log|f_1|d \log|f_2|\wedge d \log|f_3|-3\log|f_1|d\arg (f_2)\wedge d\arg(f_3).$$

On the other hand there is the canonical map $\wt {\mc L}_2\colon \mc B_2(\C)\to \R$, given by the Bloch-Wigner dilogarithm (see \cite{goncharov1995geometry}). A. Goncharov conjectured that for any algebraically closed field $k$ and any smooth projective curve $X$ over $k$, there should exist the canonical map $\mathcal H_{k(X)}\colon \Lambda^3 k(X)\to\mc B_2(k)$ such that for $k=\C$ we have $\mc P_2(X; f_1, f_2, f_3)=\wt{\mc L}_2(\mc H_{\C(X)}(f_1\wedge f_2\wedge f_3))$. The word ``canonical"  means that this map should be functorial under non-constant morphisms of curves. Moreover, motivated by the analytic properties of Chow dilogarithm, he conjectured that the map $\mc H_{k(X)}$ should additionally satisfy to the following two properties:
\begin{enumerate}
    \item It should vanish on the elements of the form $c\wedge f_2\wedge f_3, c\in k, f_2, f_3\in k(X)$,
    \item This map should give a homotopy between the map $\sum\limits_{z\in X^{(1)}}\partial_z^{(3)}$ and the zero map:
\begin{equation}
    %\label{diagram:srl_definition}
        \begin{tikzcd}[row sep=huge]
{\mathcal B_3(k(X)}&{{\mathcal B_2}(k(X))\otimes k(X)^\t} & {\Lambda^{3} k(X)^{\times}} \\
{}	&{{\mathcal B_2}(k)} & {\Lambda^{2}(k^\times).}
	\arrow["{\delta_{3}}", from=1-2, to=1-3]
	\arrow["\mc H_{k(X)}"', dashed, from=1-3, to=2-2]
		\arrow["\delta_3", from=1-1, to=1-2]
	\arrow["{-\delta_{2}}", from=2-2, to=2-3]
	\arrow["{\sum\limits_{z\in X^{(1)}}\ts_{z}}", from=1-3, to=2-3]
	\arrow["{\sum\limits_{z\in X^{(1)}}\ts_{z}}"', from=1-2, to=2-2]
\end{tikzcd}
    \end{equation}
\end{enumerate}

In this paper we assume that the field $k$ is an algebraically closed field of characteristic zero. In this case we prove the above conjecture. That is for any smooth projective curve over $k$,  we construct  the map $\mathcal H_{k(X)}\colon \Lambda^3 k(X)^\t\to\mathcal B_2(k)$ such that all these maps satisfy the conditions stated above.

\begin{remark}\begin{enumerate}
\item We impose the condition on characteristic of the field $k$ only for simplicity. It seems that the results of this paper can be generalized to the case of arbitrary characteristic. Meanwhile, the condition that $k$ is algebraically closed is essential. If the field $k$ were not algebraically closed, then in the case $k\subsetneq k(z)$ there would be no natural morphism of complexes $\Gamma(F,n)\to \Gamma(k,n-1)[-1]$. The reason is that while there is the natural map $\Gamma(F,n)\to \Gamma(k(z),n-1)[-1]$, the push-forward map $\Gamma(k(z),n-1)[-1]\to \Gamma(k,n-1)[-1]$ cannot be defined on the level of complexes. 

    \item Let $k_0$ be some subfield of $k$. It can be deduced from our main result, that if the curve $X$ together with three functions $f_1,f_2,f_3$ are defined over $k_0$ then the element $\mathcal H_{k(X)}(f_1\wedge f_2\wedge f_3)$ lies in the invariants $\mathcal B_2(k)^{Gal(k/k_0)}$. However, it seems that in general the group $\mathcal B_2(k)^{Gal(k/k_0)}$ is strictly bigger than $\mathcal B_2(k_0)$.
    
    \item By theorem of A. Suslin (Corollary 5.7 from \cite{suslin1991k3}) when the field $k$ is algebraically closed, the group $\mc B_2(k)$ is uniquely divisible. So it seems that the restriction that we work only $\Q$-linearly is not essential.
\end{enumerate}

\end{remark}

\subsection{Definitions}
\label{sec:intro:def}
We recall that everywhere we work $\mathbb Q$-linearly. 
Let $F$ be an arbitrary field. We repeat the definition of the complex $\Gamma(F, n)$ for convenience.

\begin{definition}
\label{def:pol_complexes}
Define the complex $\Gamma(F, n)$ as follows:
$$\Gamma(F,n)\colon \mathcal B_n(F)\xrightarrow{\delta_n} \mathcal B_{n-1}(F)\otimes F^\times\xrightarrow{\delta_n}\dots\xrightarrow{\delta_n}\mathcal B_2(F)\otimes \Lambda^{n-2}F^\times\xrightarrow{\delta_n}\Lambda^n F^\times.$$

This complex is concentrated in degrees $[1,n]$. The group $\mathcal B_n(F)$ is the quotient of the free abelian group generated by symbols $\{x\}_n, x\in \mathbb P^1(F)$ by some explicitly defined subgroup $\mathcal R_n(F)$ (see \cite{gon94}).  The differential is defined as follows: $\delta_n(\{x\}_k\otimes y_{k+1}\wedge \dots \wedge y_n)=\{x\}_{k-1}\otimes x\wedge y_{k+1}\wedge \dots \wedge y_n$ for $k>2$ and $\delta_n(\{x\}_2\otimes y_3\wedge\dots y_n)=x\wedge (1-x)\wedge y_3\wedge \dots \wedge y_n$.
\end{definition}

\begin{remark}\begin{enumerate}
    \item It is not known whether the definitions of the group $\mathcal R_n(F)$ from \cite{goncharov1995geometry} and \cite{gon94} are equivalent. While it is believed to be the case, this statement relies on the so-called Suslin rigidity conjecture. In this paper we use the later definition, that is the definition from \cite{gon94}.
    \item Everywhere in this paper we can replace the complex $\Gamma(F,n)$ with its canonical truncation $\tau_{\geq n-1}\Gamma(F,n)$. Therefore, only the definition of the group $\mathcal R_2(F)$ is relevant for us. As it was noted in Section 4.2 of \cite{gon94} this group  is generated by the following elements:

$$\sum\limits_{i=1}^5(-1)^i\{c.r.(x_1,\dots, \widehat x_i,\dots, x_5)\}_2, \{0\}_2, \{1\}_2, \{\infty\}_2.$$
  In this formula $x_i$ are five different points on $\mathbb P^1$ and $c.r.(\cdot)$ is the cross ratio.
\end{enumerate}

\end{remark}

Let $F$ be an arbitrary field. We recall that the $n$-th Milnor $K$-theory $K_n^M(F)$ of the field $F$ is defined as the  quotient of the vector space $\Lambda^n F^\t$ by the elements of the form $a_1\wedge (1-a_1)\wedge a_3\wedge\dots\wedge a_n,$, where $a_i\in F^\t$ and $a_1\ne 1$. We have the canonical identification $H^n(\Gamma(F,n))\cong K_n^M(F)$. If $j\colon F_1\emb F_2$ is an embedding of fields, denote by $j_*\colon K_n^M(F_1)\to K_n^M(F_2)$ the natural map given by the formula $j_*(a_1\wedge \dots\wedge a_n)=j(a_1)\wedge\dots \wedge j(a_n)$.  Bass and Tate \cite{bass1973milnor} constructed the norm map $\widehat N_{F_2/F_1} \colon K_n^M(F_2)\to K_n^M(F_1)$ which a priori depends on the choice of generators of $F_2$ over $F_1$. A. Suslin \cite{suslin1979reciprocity} proved that the norm map $\widehat N_{F_2/F_1}$ is independent of the choice of generators and is determined only by the embedding $j$.

Let $(F,\nu)$ be a discrete valuation field. Denote $\mc O_\nu=\{x\in F|\nu(x)\geq 0\}, m_\nu=\{x\in F|\nu(x)>0\}$ and $\overline F_\nu =\mc O_\nu/m_\nu$. We recall that an element $a\in F^\t$ is called \emph{a uniformiser} if $\nu(a)=1$ and \emph{a unit} if $\nu(a)=0$. For $u\in \mc O_\nu$ denote by $\overline u$ its residue class in $\oF_\nu$.

The proof of the following proposition can be found in \cite{goncharov1995geometry}:

\begin{proposition}
\label{prop:tame_symbol}
Let $(F,\nu)$ be a discrete valuation field and $n\geq 3$. There is a unique morphism of complexes $\ts[n]_\nu\colon \Gamma(F,n)\to \Gamma(\ol F_\nu, n-1)[-1]$ satisfying the following conditions:

\begin{enumerate}
    %\item For any units $u_1,\dots u_n$ we have $\ts[n]_\nu(u_1\wedge\dots \wedge u_n)=0$.
   \item For any uniformiser $\pi$ and units $u_2,\dots u_n\in F$ we have $\ts[n]_\nu(\pi\wedge u_2\wedge\dots \wedge u_n)=\overline {u_2}\wedge \dots \wedge \overline{u_n}$.
    \item For any $a\in F\bs \{0,1\}$ with $\nu(a)\ne 0$, an integer $k$ satisfying $2\leq k\leq n$ and any $b\in \Lambda^{n-k}F^\t$ we have $\ts[n]_\nu(\{a\}_k\otimes b)=0$.
    \item For any unit $u$, an integer $k$ satisfying $2\leq k\leq n$ and $b\in \Lambda^{n-k}F^\t$ we have $\ts[n]_\nu(\{u\}_k\otimes b)=-\{\overline u\}_k\otimes \ts[n-k]_\nu(b)$.
\end{enumerate}
\end{proposition}

We will call the map $\ts[n]_\nu$ from the previous proposition \emph{the tame symbol map}. The proof of this proposition can be found in \cite[Section 14]{goncharov1995geometry}. 

%\begin{remark}
%While in general the proof from \cite[Section 14]{goncharov1995geometry} is correct, the corresponding map is a morphism of complexes only up to a sign. So we need to introduce the additional minuses in items $(ii)$ and $(iii)$ of Proposition \ref{prop:tame_symbol} to make the map $\ts[n]_\nu$ a genuine morphism of complexes.
%\end{remark}

We will need the following lemma which easily follows  from the definition of the tame-symbol:

\begin{lemma}
\label{lemma:leibniz_rule_tame_symbol}
Let $(F, \nu)$ be a discrete valuation field. Let $k, n$ be two natural numbers satisfying the condition $k<n$. Let $a_1,\dots, a_n\in F^\t$ such that $\nu(a_{k+1}), \dots, \nu(a_n)=0$. Then the following formula holds:
$$\ts[n]_\nu(a_1\wedge\dots\wedge a_n)=\ts[k]_\nu(a_1\wedge\dots\wedge a_k)\wedge \overline{a_{k+1}}\wedge \dots\wedge \overline{a_n}.$$
\end{lemma}

When $D$ is an irreducible divisor on a smooth variety $X$, we denote by $\nu_D$ the corresponding discrete valuation of the field $k(X)$. For any field $F$ denote by $\nu_{\infty, F}$ the discrete valuation of $F(t)$ given by the point $\infty\in\mathbb P^1(F)$.

We recall that we have fixed some algebraically closed field $k$ of characteristic zero. Denote by $\F_d$ the category of finitely generated extensions of $k$ of transcendence degree $d$. Any morphism in this category is a finite extension. For $F\in \F_d$, denote by $\dval(F)$ the set of discrete valuations given by an irreducible Cartier divisor on some  birational model of $F$. When $F\in\F_1$ this set is equal to the set of all discrete valuations that are trivial on $k$. In this case, we denote this set simply by $\val(F)$. If $X$ is an algebraic variety together with the isomorphism $k(X)\to F$ we denote by $\dval(F)_X\subset \dval(F)$ the subset of divisorial valuations coming from divisors on $X$.

Let $j\colon K\emb F$ be an extension from $\F_d$ and $\nu\in \dval(K)$. Denote by $\ext(\nu, F)$ the set of extensions of the valuation $\nu$ to $F$. Let $\nu'\in \ext(\nu, F)$. Denote by $j_{\nn}$ the natural embedding $\oK_\nu\emb \oF_{\nu'}$. \emph{The inertia degree} $f_{\nn}$ is defined as $\deg j_{\nn}$. \emph{The ramification index} $e_{\nn}$ is defined by the formula $\pi_K=u\pi_F^{e_{\nn}}$, where $\pi_K, \pi_F$ are uniformisers of $K, F$ and $u$ is some unit. By \cite[Chapter II, \S 8]{neukirch2013algebraic} the set $\ext(\nu, F)$ is finite and, moreover, the following formula holds:
\begin{equation}
    \label{formula:degree_as_sum}
    \sum\limits_{\nu'\in \ext(\nu, F)}e_{\nn}f_{\nn}=[F:K].
\end{equation}

By the Theorem of O. Zariski \cite[Chapter VI, \S 14, Theorem 31]{zariski2013commutative} a discrete valuation on $F$ is divisorial if and only if the corresponding residue field is finitely generated over $k$ and has transcendence degree $1$. This implies that for any $\nu\in \dval(K)$, we have $\ext(\nu, F)\subset \dval(F)$. 

 For any $n\geq 0$ there is the natural map $j_*\colon \Lambda^n K^\t\to \Lambda^n F^\t$ given by the formula $j_*(a)=a$. It is easy to see that for any $\nu'\in \ext(\nu, F)$ the following formula holds:
\begin{equation}\label{formula:functoriality_of_tame_symbol}
    \ts[n]_{\nu'}j_*(a)=e_{\nn}\cdot (j_{\nn})_ *(\ts[n]_\nu(a)).
\end{equation} 

\subsection{Lifted reciprocity maps}
\label{sec:lifted_reciprocity_maps}
We recall that we work $\Q$-linearly.
%The following definition is a reformulation of Definition \ref{def:lifted_rec_map_summary}.
\begin{definition}
\label{def:SRL}

Let $F\in \F_1$. \emph{A lifted reciprocity map} on the field $F$ is a $\Q$-linear map $h\colon \Lambda^{3} F^\times\to {\P}(k)$ satisfying the following conditions:

\begin{enumerate}
    \item The following diagram is commutative:
    \begin{equation}
    \label{diagram:srl_definition}
        \begin{tikzcd}[row sep=huge]
{\mathcal B_3(F)}&{{\P}(F)\otimes F^\t} & {\Lambda^{3} F^{\times}} \\
{}	&{{\P}(k)} & {\Lambda^{2}(k^\times).}
	\arrow["{\delta_{3}}", from=1-1, to=1-2]
	\arrow["{\delta_{3}}", from=1-2, to=1-3]
	\arrow["h"', dashed, from=1-3, to=2-2]
	\arrow["{-\delta_{2}}", from=2-2, to=2-3]
	\arrow["{\sum\limits_{\nu\in \val(F)}\ts_{\nu}}", from=1-3, to=2-3]
	\arrow["{\sum\limits_{\nu\in \val(F)}\ts_{\nu}}"', from=1-2, to=2-2]
\end{tikzcd}
    \end{equation}
    \item The map $h$ vanishes on the image of the multiplication map $ \Lambda^2 F^\t\otimes  k^\t\to \Lambda^3 F^\t$.
     
\end{enumerate}

\end{definition}

\begin{remark}
\label{rem:SRL_affine}
The set of all lifted reciprocity maps has a structure of affine space over $\Q$ as any set of homotopies. 
\end{remark}

Denote by $\Set$ the category of sets. 
Define a contravariant functor $$\SRL\colon \F_1\to \Set$$ as follows. For any $F\in \F_1$ the set $\SRL(F)$ is equal to the set of all lifted reciprocity maps on $F$. If $j\colon K\emb F$ then $\SRL(j)(h_{F})$ is defined by the formula $h_K(a):=\dfrac 1{\deg j}h_{F}(j_*(a))$. It is not difficult to show that the assignment preserves identities. We will present the detailed proof that $\SRL$ is indeed a functor in Section \ref{sec:prel_results:SRL}.

\subsection{Main results}
\label{sec:main_results}

The following theorem is a solution of Conjecture 6.2 from \cite{goncharov2005polylogarithms}.

\begin{theorem}
\label{th:limit_of_functor}
For any field $F\in\F_1$ one can choose a lifted reciprocity map $\mathcal H_F$ on the field $F$ such that for any embedding $j\colon F_1\to F_2$ we have $\SRL(j)(\mathcal H_{F_2})=\mathcal H_{F_1}$. Such a collection of lifted reciprocity maps is unique.
\end{theorem}

\begin{remark}
One of the main results from \cite{rudenko_2021} states  that for any field $F\in \F_1$ there is a map $\Lambda^3 F^\t\to\P(k)$ satisfying the first condition of Definition \ref{def:SRL}. However, it is not clear why this map can be chosen functorial. The functoriality is our new result.

We remark that even the proof of existence of a homotopy is  simpler because it does not rely on complicated lemmas 5.2 - 5.7 from \cite{rudenko_2021}.
\end{remark}

\begin{remark}
In \cite[Section 6]{goncharov2005polylogarithms}, A. Goncharov proved that for any elliptic curve $E$ over $k$ there is a  lifted reciprocity map on the field $k(E)$. From the proof of Theorem \ref{th:limit_of_functor} it is not difficult to show that his map coincides with ours. Therefore, Theorem \ref{th:limit_of_functor} generalizes A. Goncharov's construction to curves of arbitrary genus.
\end{remark}

\subsubsection{Chow dilogarithm}

The definition of Chow dilogarithm can be found in Section 6 of \cite{goncharov2005polylogarithms}.  This function associates to any smooth projective curve $X$ over $\C$ and three non-zero rational functions $f_1,f_2,f_3$ on $X$ the value $\CD(X; f_1, f_2, f_3)\in \R$. The remark after Conjecture 6.2 in loc. cit. implies that Theorem \ref{th:limit_of_functor} has the following corollary:

\begin{corollary}
\label{cor:Chow_dilogarithm}
For any smooth projective curve $X$ over $\C$ and three non-zero rational functions $f_1,f_2,f_3$ on $X$ the following formula holds:
$$\CD(X; f_1, f_2, f_3)=-\wt {\mc L}_2(\mc H_{\C(X)}(f_1\wedge f_2\wedge f_3)).$$
Here $\wt {\mc L}_2\colon {\P}(\C)\to \R$ is a map given on the generators $\{x\}_2$ by the formula $$\wt{\mc L}_2(\{x\}_2)=\mc L_2(x),$$ where $\mc L_2$ is Bloch-Wigner dilogarithm.
\end{corollary}

\begin{remark}\begin{enumerate}
    \item The sign comes from the fact that we use a little bit different definition of the map $\ts_\nu$.
    \item This statement is similar to Corollary 1.5 from \cite{rudenko_2021}. However the proof from  loc. cit. is not correct: it relies on a remark after Conjecture 6.2 from \cite{goncharov2005polylogarithms}, which uses the functorial property. So Corollary \ref{cor:Chow_dilogarithm} is new. 
\end{enumerate}

\end{remark}

\subsubsection{Two-dimensional reciprocity law} From the proof of Theorem \ref{th:limit_of_functor} we get the following corollary:

\begin{corollary}
\label{cor:reciprocity_map_for_surfaces}
Let $L\in\F_2$. For any $b\in \Lambda^{4}L^\t$ and all but finitely many $\nu\in \dval(L)$ we have $\mc H_{\oL_\nu}\ts[4]_\nu(b)=0$. Moreover, the following sum is equal to zero:
\begin{equation}
\label{formula:rec_map_surfaces}
   \sum\limits_{\nu\in \dval(L)}\mc H_{\oL_\nu}\ts[4]_\nu(b)=0. 
\end{equation}

\end{corollary}

Applying $\wt{\mc L}_2$ to both sides of (\ref{formula:rec_map_surfaces}) and using Corollary \ref{cor:Chow_dilogarithm}, we recover the functional equation for Chow dilogarithm proved by A. Goncharov in \cite[Section 1.4]{goncharov2005polylogarithms}, see also \cite{gil2015simplicial}. Actually Corollary \ref{cor:reciprocity_map_for_surfaces} was our motivation behind the construction of the map $\mc H$.

%By Corollary \ref{cor:Chow_dilogarithm}, if we apply Bloch-Wigner dilogarithm to Corollary  \ref{cor:reciprocity_map_for_surfaces} we recover a functional equation for Chow dilogarithm  proved in \cite[Section 3.2]{goncharov2005polylogarithms} (see also \cite{gil2015simplicial}). If we express Chow dilogarithm through Bloch-Wigner dilogarithm using Corollary \ref{cor:Chow_dilogarithm}, the relation in loc. cit. gives some relation for Bloch-Wigner dilogarithm. Corollary \ref{cor:reciprocity_map_for_surfaces} implies that this relation is admissible (see \cite{kirillov1995dilogarithm}).

\subsection{The outline of the paper}
In Section \ref{sec:prel_results:SRL} we will show that there is the unique lifted reciprocity map on the field of rational functions $k(t)$. Denote it by $\mathcal H_{k(t)}$. Let $F\in\F_1$ be any field. Choose an embedding $j\colon k(t)\emb F$. To define $\mathcal H_F$, we extend a lifted reciprocity map $\mathcal H_{k(t)}$ from the field $k(t)$ to the field $F$. For this we solve the more general problem: for any finite extension $j'\colon F_1\emb F_2$ in $\F_1$ we construct the canonical map $N_{F_2/F_1}\colon \SRL(F_1)\to \SRL(F_2)$. More precisely we will prove the following theorem:

\begin{theorem}
\label{th:norm_map}For any embedding of fields $j\colon F_1\emb F_2$ one can define the canonical map $N_{F_2/F_1}\colon \SRL(F_1)\to \SRL(F_2)$ satisfying the following properties:
\begin{enumerate}
\item $\SRL(j)\circ N_{F_2/F_1}=id$. \item If $F_1\subset F_2\subset F_3$ is a tower of extension from $\F_1$ then $N_{F_3/F_1}=N_{F_3/F_2}\circ N_{F_2/F_1}$.
\end{enumerate}
\end{theorem}

Item $(i)$ shows that $N_{F_2/F_1}$ is indeed an extension, while item $(ii)$ shows that this extension is functorial.

Sections 2 and 3 are devoted to the proof of this theorem. The details will be given below. Now the proof of Theorem \ref{th:limit_of_functor} is easy: the existence follows from Theorem \ref{th:norm_map} together with the fact that the element $N_{F/k(t)}(\mathcal H_{k(t)})$ does not depend on the embedding $j\colon k(t)\emb F$. The uniqueness follows from standard arguments. We will present this proof in Section \ref{sec:proof_of_main_result}. The proof of Corollary \ref{cor:reciprocity_map_for_surfaces} will be given in Section \ref{sec:proof_of_reciprocity_law_for_surfaces}.

Let us outline the proof of Theorem \ref{th:norm_map}. The proof of this theorem is in many respect similar to the construction of the norm map on Milnor $K$-theory. (See \cite{bass1973milnor,suslin1979reciprocity,MILNOR1969/70, kato1980generalization}). That is the reason why we denote it by the letter $N$. (Note that compared to the norm map in Milnor $K$-theory, in our case, the norm map is directed in the opposite direction. The reason for this is that while Milnor $K$-theory gives a covariant functor, the functor $\SRL$ is contravariant.) In Section \ref{sec:norm_map:theta} for any field $F\in\F_1$ and any $\nu\in \dval(F(t))$ we will construct the map $\N_\nu\colon \SRL(F)\to \SRL(\overline {F(t)}_\nu)$. Using this map, for any extension $F_1\emb F_2$ with a generator $a$ we will define the norm map $N_{F_2/F_1, a}\colon \SRL(F_1)\to \SRL(F_2)$ (see Definition \ref{def:norm_map}). Using ideas from \cite{suslin1979reciprocity} we will show that this map does not depend on $a$ and will have finished the proof of Theorem \ref{th:norm_map}. This will be done in Section \ref{sec:norm_map:extensions_of_scalars} and Section \ref{sec:norm_map:finall_proof}.

The most non-trivial part of this paper is the construction of the map $\N_\nu$. Let us give the outline of this construction. Let $F\in \F_1$. It is useful to divide the discrete valuations of the field $F(t)$ into two classes, namely  the general valuations and the special ones (see Definition \ref{def:types_of_valuations}). For the special valuations the definition of the map $\N_\nu$ is straightforward. To reduce the definition of the map $\N_\nu$ when $\nu$ is general to the previous case, we use the notion of \emph{the lift}. Let $\nu\in \dval(F(t))$ be a general valuation and $n,j\in\mathbb N$. A lift of an element $a\in \Gamma(\overline{F(t)}_\nu, n)_j$ is an element $b\in \Gamma(F(t), n+1)_{j+1}$, such that the tame-symbol $\ts[n+1]_\nu(b)$ is equal to $a$ and the tame-symbol of $b$ at any other general valuation vanishes. The set of all lifts of the element $a$ is denoted by $\mathcal L(a)$. In Section \ref{sec:prel_results:lift} we will show that in the case $n=3, j\in\{2,3\}$, for any $a\in \Gamma(\ol {F(t)}_\nu,n)_j$ the set $\mathcal L(a)$ is non-empty. Now, when $\nu\in \dval(F(t))$ is a general valuation, $h\in \SRL(F)$ and $a\in \Lambda^3\overline{F(t)}_\nu$, we can choose some lift $b\in \mathcal L(a)$ and define the element $\N_\nu(h)(a)$ by the following formula:
$$\mathcal N_\nu(h)(a)=-\sum\limits_{\mu\in\dval(F(t))_{sp}}\mathcal N_{\mu}(h)(\ts[4]_\mu(b)).$$

Here $\dval(F(t))_{sp}$ denotes the set of all special valuations. In this formula the lifted reciprocity maps $\N_{\nu}(h)$ are already defined because $\mu$ are special. It remains to show that this expression does not depend on the choice of $b$ and for fixed $h$ gives a lifted reciprocity map on the field $\overline{F(t)}_\nu$. This can be done using the properties of the lift established in Section \ref{sec:prel_results:lift} and some version of the Parshin reciprocity law which will be proved in Section \ref{sec:prel_results:Parshin}.

\subsection{Conventions}
\label{subsec:conventions}

If $C$ is a chain complex denote by $C_d$ the elements lying in degree $d$. The symbol $\delta_n$ means the differential in the polylogarithmic complex $\Gamma(F,n)$. Although it depends on the field $F$ we will omit the corresponding sign from the notation. In the same way, when $(F,\nu)$ is a discrete valuation field we denote by $\ts[n]_\nu$ the tame-symbol map $\Gamma(F,n)\to \Gamma(\oF_\nu,n-1)[-1]$.

\section{Preliminary results}
\label{sec:prel_results}

\subsection{Lifted reciprocity maps}
\label{sec:prel_results:SRL}
\begin{proposition}\label{prop:SRL_is_functor}
$\SRL$ is indeed a functor.
\end{proposition}
\begin{proof}
If $j_1, j_2$ are some embeddings from $\F_1$ then the formula $$\SRL(j_2\circ j_1)=\SRL(j_1)\circ \SRL(j_2)$$ follows from the fact that the ramification index is multiplicative. So it is enough to show that for any embedding $j\colon K\emb F$  and $h_F\in \SRL(F)$ the map 
$$h_K:=\SRL(j)(h_F)=\dfrac 1{\deg j}h_F\circ j_*\colon \Lambda^3 K^\t\to\P(k)$$ is a lifted reciprocity map on $K$. 

The statement that $h_K$ is zero on the image of the map $\Lambda^2 K^\t\otimes k^\t\to \Lambda^3 K^\t$  follows
from the corresponding statement for $h_F$. Let us prove that diagram (\ref{diagram:srl_definition}) is commutative. 

For any $\nu\in \val(K)$ and any $\nu'\in \ext(\nu, F)$ we have $f_{\nn}=1$. Therefore, formula (\ref{formula:degree_as_sum}) becomes $\sum\limits_{\nu'\in \ext(\nu, F)}e_{\nn}=[F:K]$. Since in our case $\overline K_{\nu}\cong \overline F_{\nu'}\cong k$, the formula (\ref{formula:functoriality_of_tame_symbol}) takes the form $e_{\nn}\ts_\nu(a)=\ts_{\nu'}j_*(a)$.

For any $a\in \Lambda^3 K^\t$, we have:

\begin{equation*}
\begin{split}
    \delta_2(h_K(a))=\dfrac 1{[F:K] }\delta_2(h_F(j_*(a)))=\dfrac 1{[F:K]}\sum\limits_{\nu'\in \val(F)}\ts_{\nu'}j_*(a)=\\
         =\dfrac 1{[F:K]}\sum\limits_{\nu\in \val(K)}\sum\limits_{\nu'\in \ext(\nu, F)}\ts_{\nu'}j_*(a)=\\
         =\dfrac 1{[F:K]}\sum\limits_{\nu\in \val(K)}\sum\limits_{\nu'\in \ext(\nu, F)}e_{\nn}\ts_\nu(a)=
     \sum\limits_{\nu\in\val(K)}\ts_\nu(a).
    \end{split}
\end{equation*}
Here in the fourth equality we have used the formula $\ts_{\nu'}(j_*(a))=e_{\nn}\ts_\nu(a)$ and in the last formula we have used the formula $\sum\limits_{\nu'\in \ext(\nu, F)}e_{\nn}=[F:K]$. So the lower right triangle is commutative. The commutativity of the upper left triangle is similar. 
\end{proof}

\begin{proposition}
\label{prop:SRL_P1}
On the field $k(t)$ there is the unique lifted reciprocity map. We will denote it by $\mc H_{k(t)}$
\end{proposition}
\begin{proof}Elementary calculation shows that the group $\Lambda^3 k(t)^\t$ is generated by the image of the multiplication map $ k(t)^\t\otimes \Lambda^2 k^\t\to \Lambda^3k(t)^\times$ and by the image of $\delta_3$. Uniqueness follows from this statement.
Existence was proved in  \cite[Theorem 6.5]{goncharov1995geometry}. Let us give two remarks:
\begin{enumerate}
\item Because we use another sign convention in the definition of $\ts_\nu$ (see Proposition \ref{prop:tame_symbol}) we need to multiply the map $h$ (see Proposition 6.6) from \cite{goncharov2005polylogarithms} by $-1$.
    \item Although the proof of Proposition 6.6 from \cite{goncharov1995geometry} uses  rigidity argument, this proposition can be easily deduced from \cite{dupont1982generation}, where it was proved that ${\P}(k(t))$ is generated by elements of the form $\left\{at+b\right\}_2, a,b\in k$.
\end{enumerate}

\end{proof}

\subsection{The construction of the lift}
\label{sec:prel_results:lift}

\begin{definition}
\label{def:types_of_valuations}
Let $F\in \F_1$. A valuation $\nu\in \dval(F(t))$ is called \emph{general} if it corresponds to some irreducible polynomial over $F$. The set of general valuations are in bijection with  the set of all closed points on the affine line over $F$, which we denote by $\mb A_{F,(0)}^1$.  A valuation is called \emph{special} if it is not general. Denote the set of general (resp. special) valuations by $\dval(F(t))_{gen}$ (resp. $\dval(F(t))_{sp}$). 
\end{definition}

\begin{remark}
\label{rem:types_of_valuations}
Let us realize $F$ as a field of fractions on some smooth projective curve $X$ over $k$. Set $S=X\t \mathbb P^1$. It can be checked that a valuation $\nu\in \dval(k(S))$ is special in the following two cases:

\begin{enumerate}
    \item There is a birational morphism $p\colon \wt S\to S$, and the valuation $\nu$ corresponds to some irreducible divisor $D\subset \wt S$ contracted under $p$.
    \item The valuation $\nu$ corresponds to some of the divisors $X\t\{\infty\}, \{a\}\t\mathbb P^1, a\in X$.
\end{enumerate}

Otherwise, the valuation $\nu$ is general.
It follows from this description that if $\nu$ is a special valuation different from $X\times \{\infty\}$, then the residue field $\overline{F(t)}_\nu$ is isomorphic to $k(t)$.
\end{remark}

\begin{definition}
\label{def:lift}
Let $F\in \F_1$, $j, n\in \mathbb N$, $\nu\in \dval(F(t))_{gen}$ and $a\in \Gamma(\ol {F(t)}_\nu, n)_j$. \emph{A lift} of the element $a$ is an element $b\in \Gamma(F(t), n+1)_{j+1}$ satisfying the following two properties:
\begin{enumerate}
    \item $\ts[n+1]_\nu(b)=a$ and
    \item for any general valuation $\nu'\in\dval(F(t))_{gen}$ different from $\nu$, we have $\ts[n+1]_{\nu'}(b)=0$.
\end{enumerate}

The set of all lifts of the element $a$ is denoted by $\mathcal L(a)$.
\end{definition}

The main result of this section is the following statement:

\begin{theorem}
\label{th:main_exact_sequence}
Let $F\in \F_1$ and $j\in\{2,3\}$. For any $\nu\in \dval(F(t))_{gen}$ and $a\in \Gamma(\ol{F(t)}_\nu, 3)_j$, the following statements hold:
\begin{enumerate}
    \item  The set $\mathcal L(a)$ is non-empty.
    \item Let us assume that $j=3$. For any $b_1,b_2\in \mathcal L(a)$, the element $b_1-b_2$ can be represented in the form $a_1+\delta_4(a_2)$, where $a_1\in \Lambda^4 F^\t$ and $a_2\in \Gamma(F(t),4)_3$ such that for any $\nu\in\dval(F(t))_{gen}$ the element $\ts[4]_\nu(a_2)$ lies in the image of the map $\delta_3\colon \Gamma(F(t),3)_1\to \Gamma(F(t),3)_2$.
\end{enumerate}
\end{theorem}

\begin{remark}
Item $(i)$ shows that any element has some lift, while item $(ii)$ shows that in the case $j=3$ up to some specific elements the lift is unique. 
\end{remark}

The proof of this theorem was inspired by exact sequences in \cite{rudenko_2021}.  We need two lemmas.

\begin{lemma}
\label{lemma:Res_is_surjective}
Let $m\geq 3$ be an integer. The following map is surjective in degrees $m-1, m$:
$$\Gamma(F(t),m)\xrightarrow{\left(\ts[m]_{\nu_p}\right)} \bigoplus\limits_{p\in \mb A_{F,(0)}^1}\Gamma(F(p),m-1)[-1].$$
In this formula $\nu_p$ is the valuation corresponding to the point $p$ and $F(p)\cong \ol{F(t)}_{\nu_p}$ is the residue field of this point.
\end{lemma}

The proof of this lemma is completely similar to the proof of surjectivity in the exact sequence of Bass and Tate from \cite{bass1973milnor} describing the Milnor $K$-theory of rational function field in one variable.

\begin{proof}
For a point $p\in\mathbb A^1_{F,(0)}$ denote by $f_p\in F[t]$ the corresponding monic irreducible polynomial. Define an increasing filtration $\mathcal F$ on the complex $\Gamma(F(t),m)$ as follows: the subspace $\mathcal F_d(\Gamma(F(t),m))$ is equal to the set of elements lying in the kernels of all the maps $\ts[m]_{\nu_p}$ with $\deg p>d$. It is enough to prove that for any $d\geq 0$ the following map is surjective:
$$\mathcal F_d(\Gamma(F(t),m))\to \bigoplus\limits_{\substack{p\in \mb A_{F,(0)}^1\\\deg p\leq d}}\Gamma(F(p),m-1)[-1].$$

The proof is by induction on $d$. The case $d=-1$ is trivial. Let us prove the inductive step. It is enough to show that for any $a\in \Gamma(F(p),m-1)_{j-1}, j\in\{m-1,m\}$ there is an element $\wt a\in \Gamma(F(t),m)_{j}$ with the following properties:
\begin{enumerate}
    \item for any $p'\ne p$ with $\deg p'\geq\deg p$ we have $$\ts[m]_{\nu_{p'}}(\wt a)=0.$$
    \item We have $\ts[m]_{\nu_p}(\wt a)=a$.
\end{enumerate}

For an element $\xi\in F(p)$ there is a unique polynomial $l_p(\xi)$ of degree $<\deg p$ such that the image of $l_p(\xi)$ under the natural projection $F[t]\to F[t]/f_p\cong F(p)$ is equal to $\xi$.

 The following formulas for $\wt a$ are taken from \cite[Section 5.2]{rudenko_2021}.

\begin{description}
\item[Case $j=m-1$] Choose a representation $a=\sum\limits_{\alpha}n_{\alpha}\cdot(\{\xi_1^{\alpha}\}_2\otimes \xi_3^{\alpha}\wedge \dots \wedge \xi_{m-1}^{\alpha}) $. Define the element $\wt a$ by the formula
$$\wt a=-\sum\limits_{\alpha}n_{\alpha}\cdot (\{l_p(\xi_1^{\alpha})\}_2\otimes f_p\wedge l_p(\xi_3^{\alpha})\wedge \dots\wedge l_p(\xi_{m-1}^{\alpha})).$$
\item[Case $j=m$] Choose a representation $a=\sum\limits_{\alpha}n_{\alpha}\cdot(\xi_1^{\alpha}\wedge\dots\wedge \xi_{m-1}^{\alpha}) $. The element $\wt a$ is defined by the formula
$$\wt a=\sum\limits_{\alpha}n_{\alpha}\cdot (f_p\wedge l_p(\xi_1^{\alpha})\wedge\dots\wedge l_p(\xi_{m-1}^{\alpha})).$$
\end{description}
It is easy to see that these elements satisfy the conditions stated above.
\end{proof}

\begin{proposition}
\label{prop:Rudenko_result}
The following sequence is exact for $j=m$ and exact in the third term for $j=m-1$:
\begin{equation}
    \begin{split}
        0\to H^j(\Gamma(F,m))\to H^j(\Gamma(F(t),m))\xrightarrow{(\ts[m]_{\nu_p})}\bigoplus\limits_{p\in \mb A_{F,(0)}^1}H^{j-1}(\Gamma(F(p),m-1))\to 0.
    \end{split}
\end{equation}
\end{proposition}

\begin{proof}
We recall that for any field $K$ we have the canonical identification $H^m(\Gamma(K,m))\cong K_m^M(K)$. So in the case $j=m$ the statement of the proposition is equivalent to the following exact sequence:
$$0\to K_m^M(F)\to K_m^M(F(t))\xrightarrow{(\ts[m]_{\nu_p})}\bigoplus\limits_{p\in \mb A_{F,(0)}^1}K_{m-1}^M(F(p))\to 0.$$

This is the exact sequence of Bass and Tate from \cite{bass1973milnor} describing the Milnor $K$-theory of rational function field in one variable. The exactness in the last term for $j=m-1$ is a particular case of the main result from \cite{rudenko_2021}.
\end{proof}

\begin{proof}[The proof of Theorem \ref{th:main_exact_sequence}]

\begin{enumerate}
    \item Follows from Lemma \ref{lemma:Res_is_surjective} for $m=4$.
    \item Denote $b=b_1-b_2$. The element $b$ satisfies  $\partial_\nu^{(4)}(b)=0$ for any $\nu\in\dval(F(t))_{gen}$. Denote by $\overline{b}$ the corresponding element in $H^4(\Gamma(F(t),4))=K_4^M(F(t))$. Proposition \ref{prop:Rudenko_result} for $j=m=4$ implies that $\overline b$ lies in the image of the map $H^4(\Gamma(F,4))\to H^4(\Gamma(F(t),4))$. This implies that $b$ can be represented in the form $a_1+\delta_4(\wt{a_2})$, where $a_1\in \Lambda^4 F^\times$ and $\wt{a_2}\in \Gamma(F(t),4)_3$. It remains to show that there is $a_2\in \Gamma(F(t),4)_3$ satisfying the following two conditions:
    \begin{enumerate}
        \item $\delta_4(a_2)=\delta_4(\wt {a_2})$ and
        \item $\ts[4]_\nu(a_2)$ lies in the image of $\delta_3$.
    \end{enumerate}

    For any $\nu\in \dval(F(t))_{gen}$ we set $b_\nu=\partial_\nu^{(4)}(\wt{a_2})$. We have: $$\delta_3(b_\nu)=\delta_3(\partial_\nu^{(4)}(\wt{a_2}))=-\partial_\nu^{(4)}(\delta_4(\wt{a_2}))=-\partial_\nu^{(4)}(b-a_1)=-\partial_\nu^{(4)}(b)+\partial_\nu^{(4)}(a_1)=0.$$ So $b_\nu$ lies in the kernel of $\delta_3$ and gives the element $\ol{b_\nu}\in H^2(\Gamma(\overline {F(t)}_\nu,3))$. Consider the element $$(\ol{b_\nu})_{\nu\in \dval(F(t))_{gen}}\in \bigoplus\limits_{\nu\in \dval(F(t))_{gen}}H^2(\Gamma(\overline {F(t)}_\nu,3)).$$ Proposition \ref{prop:Rudenko_result} for $j=3, m=4$ shows that there is an element $a_3\in H^3(\Gamma(F(t),4))$ such that for any $\nu\in\dval(F(t))_{gen}$ we have $\ts[4]_\nu(a_3)=\ol{b_\nu}$. Let $\wt a_3$ be arbitrary lift of $a_3$ to $\Gamma(F(t),4)_3$  and $a_2=\wt{a_2}-\wt a_3$. By construction, for any $\nu\in \dval(F(t))_{gen}$, the element $\ts[4]_\nu(a_2)$ is zero in $H^2(\Gamma(\ol{F(t)}_\nu,3))$ and hence lies in the image of $\delta_3$. Since $\wt{a_3}$ lies in the kernel of $\delta_4$, we have $\delta_4(a_2)=\delta_4(\wt{a_2})$. So the two conditions above hold.

\end{enumerate}
\end{proof}

\subsection{Parshin reciprocity law}
\label{sec:prel_results:Parshin}
The goal of this section is to prove the following theorem:

\begin{theorem}
\label{th:Parshin_sum_point_curve}
Let $L\in\F_2$ and $j\in\{3,4\}$. For any $b\in \Gamma(L,4)_j$ and all but finitely many $\mu\in \dval(L)$ the following sum is zero:
$$\sum\limits_{\mu'\in \val(\oL_\mu)}\ts_{\mu'}\ts[4]_\mu(b)=0.$$
Moreover the following sum is zero:
$$\sum\limits_{\mu\in \dval(L)}\sum\limits_{\mu'\in \val(\oL_\mu)}\ts_{\mu'}\ts[4]_\mu(b)=0.$$
\end{theorem}

%\begin{remark}
%The assertion will remain true if we replace the group $\Gamma(L,4)_j$ with the group $\Gamma(L,m)_j$ for arbitrary $m\geq 3, 1\leq j\leq m$. In the case $m=j=3$ we will get some version of the classical Parshin reciprocity law (see \cite{parshin1975class}, \cite{osipov2011categorical}) for strictly regular elements.
%\end{remark}

Let $X$ be a smooth algebraic variety. The definition of a simple normal crossing divisor can be found in \cite{kollar2009lectures}. Denote by $X^{(1)}$ the set of all closed irreducible subsets of $X$ of codimension $1$. For a divisor $D=\sum\limits_{D\in X^{(1)}}n_D[D]$ on $X$, denote by $|D|$ \emph{its support} defined by the formula $\sum\limits_{\substack {D\in X^{(1)}\\n_D\ne 0}} [D]$. A divisor $D$ is called \emph{supported on a simple normal crossing divisor} if $|D|$ is a simple normal crossing divisor. Let $x\in X$. A divisor $D$ is called \emph{supported on a simple normal crossing divisor locally at $x$}, if the restriction of this divisor to some open affine neighborhood of the point $x$ is supported on a simple normal crossing divisor.

We have the following statement:

\begin{theorem}
\label{th:resolution_of_singularity}
Let $X$ be a variety over an algebraically closed field of characteristic zero and $D$ an effective Weil divisor on $X$. There is a birational morphism $f\colon \wt X\to X$ such that $\wt X$ is smooth and $f^{*}(D)$ is supported on a simple normal crossing divisor.
\end{theorem}

 For a rational function $f$ on $X$, denote by $(f)$ its divisor. For the definition of the complexes $\Gamma(F,n)$, see Definition \ref{def:pol_complexes}.
\begin{definition}
Let $X$ be a smooth algebraic variety and $x\in X$. An element of the vector space $\Gamma(k(X),4)_3$ (resp. $\Gamma(k(X),4)_4$) is called \emph{strictly regular at $x\in X$} if it can be represented as a linear combination of elements of the form $\{\xi_1\}_2\otimes \xi_3\wedge \xi_4$ (resp. $\xi_1\wedge \xi_2\wedge \xi_3\wedge  \xi_4$) such that all the divisors $|(\xi_1)|+|(\xi_3)|+|(\xi_4)|$ (resp. $|(\xi_1)|+|(\xi_2)|+|(\xi_3)|+|(\xi_4)|$) are supported on a simple normal crossing divisor locally at $x$. 
\end{definition}

Theorem \ref{th:resolution_of_singularity} has the following corollary:

\begin{corollary}
\label{cor:existence_of_birational_morphism}
Let $S$ be a smooth surface and $j\in \{3,4\}$. For any element $a\in \Gamma(k(S),4)_j$ there is a birational morphism $p\colon \wt S\to S$ such that the element $p^*(a)$ is strictly regular at all points.
\end{corollary}

The following lemma characterizes strictly regular elements:

\begin{lemma}
\label{lemma:characterisation_of_strictly_regular_elements}
Let $S$ be a smooth algebraic surface and $x\in S$. 

\begin{enumerate}
    \item The subgroup of strictly regular elements of $\Gamma(k(S),4)_3$ is generated by elements of the following form:
    
\begin{enumerate}
    \item $\{\pi_1^n\pi_2^m \xi_1\}\otimes \pi_1\wedge \pi_2$.
    \item $\{\pi_1^n\pi_2^m \xi_1\}\otimes \pi_1\wedge \xi_4$.
    \item $\{\pi_1^n\pi_2^m \xi_1\}\otimes \xi_3\wedge\xi_4$.
\end{enumerate}

Here all the functions $\xi_i$ take non-zero values at $x$ and $\pi_i$ is a regular system of parameters.
    
    \item The subgroup of strictly regular elements of $\Gamma(k(S),4)_4$ is generated by elements of the following form:
    
\begin{enumerate}
    \item $\pi_1\wedge\pi_2\wedge \xi_3\wedge\xi_4$.
    \item $\pi_1\wedge \xi_2\wedge \xi_3 \wedge \xi_{4}$.
    \item $\xi_1\wedge \xi_2\wedge \xi_3\wedge \xi_4$.
\end{enumerate}
The functions $\xi_i,\pi_i$ satisfy the same conditions as in item (i).
\end{enumerate}
\end{lemma}

\begin{proof}
Follows from the fact that if $\pi_1,\pi_2$ is a regular system of parameters at $x$ then any function $f\in k(S)$ can be written in the form $\pi_1^{n_1}\pi_2^{n_2}\xi$, where $n_i\in\Z$ and $\xi$ is a regular function at $x$ such that $\xi(x)\ne 0$.
\end{proof}

%\subsection{Parshin reciprocity law}

The following result is a version of the classical Parshin reciprocity law for strictly regular elements (see \cite{parshin1975class,horozov2014reciprocity, osipov2011categorical}). 

\begin{theorem}
\label{th:Parshin_reciprocity_law}
Let $S$ be a surface smooth at some point $x\in S$ and $j\in \{3,4\}$. For any strictly regular element $b$ of the group $\Gamma(k(S),4)_j$ at $x$ the following sum is equal to zero:
$$\sum\limits_{\substack{C\subset S\\ C\ni x}}\ts_{\nu_{x, C}}\ts[4]_{\nu_{C}}(b)=0. $$
Here the sum is taken over all irreducible curves $C\subset S$ containing the point $x$  and smooth at it, $\nu_C$ is the valuation corresponding to $C$ and $\nu_{x, C}$ is a valuation of the residue field $\overline{k(S)}_{\nu_C}\cong k(C)$ corresponding to $x\in C$.
\end{theorem}

\begin{proof}
It is enough to prove this theorem for any of the generators from Lemma \ref{lemma:characterisation_of_strictly_regular_elements}. We will only consider the most interesting case $(i), (a)$. We can assume that $S$ is a smooth surface, $x\in S$ and $\pi_1, \pi_2$ is a system of regular parameters at $x$. Passing to some open affine neighborhood of the point $x$, we can assume that the following conditions hold:
\begin{enumerate}
    \item the function $\xi_1$ is invertible and regular,
    \item the functions $\pi_i$ are regular,
    \item for any $i \in\{1,2\}$ the divisor of the function $\pi_i$ is equal to some irreducible curve $C_i$ passing through $x$.
\end{enumerate}
In general, if $X$ is a subvariety of an algebraic variety $Y$ and $f$ is a regular function on $Y$ we denote its restriction to $X$ by $\left.f\right|_X$. 

Let $C\subset S$ be a curve and $f_1, f_2, f_3$  regular functions on $S$. Assume that $\nu_C(f_3)=0$. It follows from Proposition \ref{prop:tame_symbol} that the tame-symbol $\ts[4]_C(\{f_1\}_2\otimes f_2\wedge f_3)$ is equal to zero if $\nu_C(f_1)\ne 0$ and is equal to $-\nu_C(f_2)\cdot(\{\left.f_1\right|_C\}_2\otimes (\left.f_3\right|_C))$ if $\nu_C(f_1)=0$.

Let $b$ be $\{\pi_1^n\pi_2^m \xi_1\}_2\otimes \pi_1\wedge \pi_2$ and $C$ be an irreducible curve on $X$. It follows from the last paragraph that  the only curves on $S$ satisfying $\ts[4]_{\nu_C}(b)\ne 0$ are $C_1$ and $C_2$. Consider the following three cases:
\begin{description}
\item[Case $n,m\ne 0$] In this case the integers $\nu_{C_i}(\pi_1^n\pi_2^m \xi_1), i=1,2$ are non-zero and so both of tame-symbols $\ts[4]_{\nu_{C_1}}(b)$ and $\ts[4]_{\nu_{C_2}}(b)$ vanish. The statement follows.
\item [Case $n\ne 0, m=0$ or $m\ne 0, n=0$] Consider, say, the first case. As in the previous item, $$\ts[4]_{\nu_{C_1}}(b)=0.$$ So it is enough to prove that $\ts_{\nu_{x,C_2}}\ts[4]_{\nu_{C_2}}(b)=0$. We have:
$$\ts[4]_{\nu_{C_2}}(b)=\{\left.(\pi_1^n\xi_1)\right|_{C_2}\}_2\otimes \left.\pi_1\right|_{C_2}.$$

Now the statement follows from the following formula: $\ord_x(\left.(\pi_1^n\xi_1)\right|_{C_2})=n\ne 0.$
\item[Case $n=m=0$] We have: $\ts[4]_{\nu_{C_1}}(b)=-\{\left.\xi_1\right|_{C_1}\}_2\otimes \left.\pi_2\right|_{C_1}.$ So $\ts[3]_{\nu_{x,C_1}}\ts[4]_{\nu_{C_1}}(b)=\{\xi_1(x)\}_2.$ Similarly, $\ts[3]_{\nu_{x,C_2}}\ts[4]_{\nu_{C_2}}(b)=-\{\xi_1(x)\}_2.$ The statement follows.
\end{description}

\end{proof}

%The previous theorem also holds in the case $q=1, j=2$. In this case we recover classical Parshin reciprocity law \cite{parshin1975class}(see also \cite{osipov2011categorical}) in the case of strictly regular elements. This explains the name of this theorem.

\begin{proof}[The proof of Theorem \ref{th:Parshin_sum_point_curve}]
Let $S$ be an algebraic surface with $k(S)\cong L$. The definition of the set $\dval(L)_S$ was given in the Section \ref{sec:intro:def}. We recall that this is the set of all discrete valuations coming from divisors on $S$. Choose $S$ in such a way that $b$  would be strictly regular at all points of $S$.  This is possible by Corollary \ref{cor:existence_of_birational_morphism}. Theorem \ref{th:Parshin_reciprocity_law} implies the following formula:
$$\sum\limits_{\mu\in \dval(L)_S}\sum\limits_{\mu'\in \val(\oL_\mu)}\ts_{\mu'}\ts[4]_\mu(b)=0.$$
It remains to prove that for any $\mu\in \dval(L)\bs \dval(L)_S$ the following sum vanishes:
$$\sum\limits_{\mu'\in \val(\oL_\mu)}\ts_{\mu'}\ts[4]_\mu(b)=0.$$
There is a birational morphism $p\colon \wt S\to S$ such that $\mu$ is given by a divisor on $\wt S$ contracted under $p$. The morphism $p$ is a sequence of blow-ups $p_m\circ\dots \circ p_1, p_i\colon S_i\to S_{i-1}, S_m=\wt S, S_0=S$. Let $D_i\subset S_i$ be the corresponding exceptional curve. Denote by $\mu_i$ the corresponding valuation. It is enough to show that for any $i$ the following formula holds:
$$\sum\limits_{\mu'\in \val(\oL_{\mu_i})}\ts_{\mu'}\ts[4]_{\mu_i}(b)=0.$$
This formula follows from Theorem \ref{th:Parshin_reciprocity_law} for the element $b$ and the surfaces $S_i$ and $S_{i-1}$. 
\end{proof}

\subsection{Lemma about finiteness}
The goal of this section is to prove the following lemma which we will need later:

\begin{lemma}
\label{lemma:finitnes_of_sum}
Let $L\in\F_2$. For any $b\in \Lambda^{4} L^\t$ and all but a finite number of $\nu\in \dval(L)$ the element $\ts[4]_\nu(b)$ belongs to the image of the multiplication map $\overline L_\nu^\t\otimes \Lambda^2 k^\t\to \Lambda^3 \overline L_\nu^\t$.
\end{lemma}

\begin{proof}
By Corollary \ref{cor:existence_of_birational_morphism}, there is a smooth proper algebraic surface $S$ such that $L=k(S)$ and $b$ is strictly regular at all points of $S$. We recall that the set $\dval(L)_S$ was defined in the Section \ref{sec:intro:def}. 

We can assume that the element $b$ has the form $f_1\wedge f_2\wedge f_3\wedge f_4$. Let $W=\sum\limits_{j=1}^4|(f_j)|$. According to Proposition \ref{prop:tame_symbol}, the tame-symbol $\ts[4]_{\nu_C}(b)$ vanishes if the curve  $C$ does not belong to divisor $W$. So for all but finitely many $\nu\in\dval(L)_S$ the tame-symbol $\ts[4]_\nu(b)$ vanishes. It remains to show that for any $\nu\in \dval(L)\backslash\dval(L)_S$, the element $\ts[4]_\nu(b)$ lies in the image of the map $\overline L_\nu^\t\otimes \Lambda^2 k^\t\to \Lambda^3 \overline L_\nu^\t$. (See also the proof of Theorem \ref{th:Parshin_sum_point_curve}).

Let $p\colon \wt S\to S$ be a birational morphism, such that $\nu$ correspond to some irreducible divisor $D$ on $\wt S$ contracted under $p$. Denote by $x\in S$ the image of $D$ under $p$. Lemma \ref{lemma:characterisation_of_strictly_regular_elements} for the surface $S$, the point $x$ and the element $b$ implies that $b$ can be represented in the form $\sum\limits_{i}m_i\cdot(f_1^{(i)}\wedge f_2^{(i)}\wedge \xi_3^{(i)}\wedge \xi_4^{(i)})$, such that the divisors of the functions $\xi_j^{(i)}$ do not contain the point $x$. This implies that $\nu_D(\xi_j^{(i)})=0$ and moreover the restrictions of the functions $\xi_j^{(i)}$, considered as functions on $\wt S$, to $D$ lie in $k$. Now the statement follows from Lemma \ref{lemma:leibniz_rule_tame_symbol} for $k=2, n=4$.
\end{proof}

\section{The norm map}
\label{sec:norm_map}
\subsection{The definition of $\N_\nu$}
\label{sec:norm_map:theta}
Let $F\in \F_1$, $\nu\in\dval(F(t))$. Denote the field $F(t)$ by $L$. The goal of this section is to construct the map $\N_\nu\colon \SRL(F)\to \SRL(\oL_\nu)$. We will do this in the following three steps:

\begin{enumerate}
    \item We will define this map when $\nu$ is a special valuation (see Definition \ref{def:types_of_valuations}).
    \item Using the construction of the lift from Section \ref{sec:prel_results:lift}, for any general valuation $\nu\in\dval(L)$, we will define the map $\N_\nu\colon \SRL(F)\to\Hom(\Lambda^3 \oL_\nu^\times, \mc B_2(k))$. (Proposition \ref{prop:sigma_well_defined}).
    \item Using Theorem \ref{th:Parshin_sum_point_curve} from Section \ref{sec:prel_results:Parshin}, we will show that for any $h\in\SRL(F)$, the map $\N_\nu(h)\colon \Lambda^3 \oL_\nu^\t\to \P(k)$ is a lifted reciprocity map on the field $\oL_\nu$. (Proposition \ref{prop:theta_is_SRL}). So, $\N_\nu$ gives a map $\SRL(F)\to \SRL(\ol{L}_\nu)$.
\end{enumerate}

We recall that the discrete valuation $\nu_{\infty, F}\in\dval(F(t))$ was defined in Section \ref{sec:intro:def}. Let $\nu$ be special. If $\nu=\nu_{\infty,F}$ then define $\N_\nu(h)=h$ (here we have used the identification of $\oL_{\nu_{\infty,F}}$ with $F$). In the other case we have $\overline{L}_\nu\simeq k(t)$ (see Remark \ref{rem:types_of_valuations}). In this case define $\N_\nu(h)$ to be the unique lifted reciprocity map from Proposition \ref{prop:SRL_P1}. We have defined $\N_\nu$ for any $\nu\in \dval(L)_{sp}$. 

Let $h\in\SRL(F)$. Define the map $H_h\colon \Lambda^{4} L^\t\to {\P}(k)$ by the following formula:
$$H_h(b)=-\sum\limits_{\mu\in \dval(L)_{sp}}\N_\mu(h)(\ts[4]_{\mu}(b)).$$ 
This sum is well-defined by Lemma \ref{lemma:finitnes_of_sum}. Recall that we defined the notion of the lift in Section \ref{sec:prel_results:lift}.

\begin{definition}
\label{def:theta_gen}
Let $\nu\in\dval(L)_{gen}$. Define the map $\mc N_\nu\colon \SRL(F)\to \Hom(\Lambda^3\oL_\nu^\t, \mathcal B_2(k))$ as follows. Let $h\in\SRL(F)$ and $a\in \Lambda^3\oL_\nu^\t$. Choose some lift $b\in\mc L(a)$ and define the element $\N_\nu(h)(a)$ by the formula $H_h(b)$. 
\end{definition}

\begin{proposition}
\label{prop:sigma_well_defined}
The previous definition is well-defined i.e. for any $b_1,b_2\in \mathcal L(a)$ we have $H_h(b_1)=H_h(b_2)$.
\end{proposition}
\begin{proof}
We need to show that the element $H_h(b_1)-H_h(b_2)=H_h(b_1-b_2)$ is equal to zero. By item $(ii)$ of Theorem \ref{th:main_exact_sequence}, it is enough to show that the map $H_h$ vanishes on the elements of the form $a_1, \delta_4(a_2)$, where $a_1\in\Lambda^4 F^\t$ and $ a_2\in\Gamma(L,4)_3$ such that for any $\nu\in\dval(L)_{gen}$ the element $\ts[4]_\nu(a_2)\in \Gamma(\oL_\nu,3)_2$ lies in the image of the map $\delta_3\colon \Gamma(\oL_\nu, 3)_1\to \Gamma(\oL_\nu, 3)_2$. 

\begin{enumerate}
    \item  Direct computation shows that for any $a_1\in \Lambda^4 F^\t$ and any $\mu\in \dval(L)$ the element $\ts[4]_\mu(a_1)$ lies in the subgroup $\Lambda^3 k^\t\subset\Lambda^{3}\oL_\mu^\t$. It follows that $H_h(a_1)=0$.
    \item  We need to show that $H_h(\delta_4(a_2))=0$. Let $\mu$ be a special valuation. By Proposition \ref{prop:tame_symbol} and the fact that the map $\N_\mu(h)$ is a lifted reciprocity map(see Definition \ref{def:SRL}) we have $\N_\mu(h)\ts[4]_\mu\delta_4(a_2)=-\N_\mu(h)\delta_3\ts[4]_\mu(a_2)=-\sum\limits_{\mu'\in \val(\oL_\mu)}\ts_{\mu'}\ts[4]_\mu(a_2)$. So by the definition of the map $H_h$ we get:
\begin{equation}
    \begin{split}
        H_h(\delta_{4}(a_2))=-\sum\limits_{\mu\in \dval(L)_{sp}}\N_\mu(h)\ts[4]_\mu\delta_{4}(a_2)=\sum\limits_{\mu\in \dval(L)_{sp}}\sum\limits_{\mu'\in \val(\oL_\mu)}\ts_{\mu'}\ts[4]_{\mu}(a_2)=\\
        = \sum\limits_{\mu\in \dval(L)}\sum\limits_{\mu'\in \val(\oL_\mu)}\ts_{\mu'}\ts[4]_{\mu}(a_2)=0.
    \end{split}
\end{equation}
Here the third equality holds because for any general valuation $\mu$ the element $\ts[4]_\mu(a_2)\in \Gamma(\oL_\mu,3)_2$ lies in the image of the map $\delta_3\colon \Gamma(\oL_\mu, 3)_1\to \Gamma(\oL_\mu, 3)_2$ and so it lies in the kernel of all the maps $\ts[3]_{\mu'}, \mu'\in\val(\oL_\mu)$. This follows from the fact that $\ts[3]_{\mu'}$ is a morphism of complexes. The fourth equality follows from Theorem \ref{th:Parshin_sum_point_curve}.
\end{enumerate}

\end{proof}

It remains to prove that for any $h\in\SRL(F)$, the map $\N_\nu(h)$ is a lifted reciprocity map on the field $\oL_\nu$. For this we need the following lemma:

\begin{lemma}\label{lemma:double_sum}
Let $j\in \{2,3\}$, $\nu\in\dval(L)_{gen}, a\in \Gamma(\oL_\nu,3)_j$ and $b\in \mc L(a)$. The following formula holds:
$$\sum\limits_{\nu'\in \val(\oL_\nu))}\ts_{\nu'}(a)=-\sum\limits_{\mu\in \dval(L)_{sp}}\sum\limits_{\mu'\in \val(\oL_{\mu})}\ts_{\mu'}\ts[4]_\mu(b).$$
\end{lemma}

This lemma allows us to reduce some statement about the field $\oL_\nu$ for $\nu\in\dval(L)_{gen}$ to the corresponding statements for the fields $\oL_\mu$ for special $\mu$. 
\begin{proof}
Theorem \ref{th:Parshin_sum_point_curve} implies the following formula:

$$\sum\limits_{\mu\in \dval(L)}\sum\limits_{\mu'\in \val(\oL_\mu)}\ts_{\mu'}\ts[4]_\mu(b)=0.$$

On the other hand:
\begin{equation*}
    \begin{split}
        \sum\limits_{\mu\in \dval(L)}\sum\limits_{\mu'\in \val(\oL_\mu)}\ts_{\mu'}\ts[4]_\mu(b)=\\
        =\sum\limits_{\nu'\in \val(\oL_\nu))}\ts_{\nu'}(a)+\sum\limits_{\mu\in \dval(L)_{sp}}\sum\limits_{\mu'\in \val(\oL_\mu)}\ts_{\mu'}\ts[4]_\mu(b).
    \end{split}
\end{equation*}

The statement of the lemma follows.
\end{proof}

\begin{proposition}\label{prop:theta_is_SRL}
Let $h\in\SRL(F)$. The map $\N_\nu(h)$ is a lifted reciprocity map.
\end{proposition}

\begin{proof}
Let us show that the following diagram is commutative:
    \begin{equation*}
        \begin{tikzcd}[row sep=huge]
{{\P}(\overline L_\nu)\otimes \overline L_\nu^\t} & {\Lambda^{3} \overline L_\nu^{\times}} \\
	{{\P}(k)} & {\Lambda^{2}(k^\times).}
	\arrow["{\delta_{3}}", from=1-1, to=1-2]
	\arrow["\N_\nu(h)"', dashed, from=1-2, to=2-1]
	\arrow["{-\delta_{2}}", from=2-1, to=2-2]
	\arrow["{\sum\limits_{\nu'\in \val(\overline L_\nu)}\ts_{\nu'}}", from=1-2, to=2-2]
	\arrow["{\sum\limits_{\nu'\in \val(\overline L_\nu)}\ts_{\nu'}}"', from=1-1, to=2-1]
\end{tikzcd}
    \end{equation*}

\begin{description}

\item[The lower right triangle]
Let $a\in \Lambda^{3} \oL_\nu^\t$. Choose some $b\in\mc L(a)$. We have:
%\begin{equation*}
    \begin{align*}
    &\sum\limits_{\nu'\in \dval(\oL_\nu))}\ts_{\nu'}(a)=-\sum\limits_{\mu\in \dval(L)_{sp}}\sum\limits_{\mu'\in \val(\oL_\mu)}\ts_{\mu'}\ts[4]_\mu(b)\quad\text{(by Lemma \ref{lemma:double_sum})}    \\
        &=\sum\limits_{\mu\in \dval(L)_{sp}} \delta_2\N_{\mu}(h)\ts[4]_\mu(b)\quad\text{(because $\N_\mu(h)$ is a lifted reciprocity map)}\\
        &=
        -\delta_2H_h(b)\quad\text{(follows from the definition of the map $H_h$)}\\
        &=-\delta_2\N_\nu(h)(a).
        \end{align*}
%\end{equation*}

\item[The upper left triangle] Let $a\in \Gamma(\oL_\nu,3)_2$. Choose $b\in \mc L(a)$. We have:

%\begin{equation*}
    \begin{align*}
    &\sum\limits_{\nu'\in \val(\oL_\nu))}\ts_{\nu'}(a)=-\sum\limits_{\mu\in \dval(L)_{sp}}\sum\limits_{\mu'\in \val(\oL_{\mu})}\ts_{\mu'}\ts[4]_\mu(b)\quad\text{(by Lemma \ref{lemma:double_sum})}\\
       & =-\sum\limits_{\mu\in \dval(L)_{sp}} \N_\mu(h)\delta_{3}\ts[4]_\mu(b)\quad\text{(because $\N_\mu(h)$ is a lifted reciprocity map)}\\
       & =
        \sum\limits_{\mu\in \dval(L)_{sp}} \N_\mu(h)\ts[4]_\mu\delta_{4}(b)\quad\text{(because $\ts[4]_\mu$ is a morphism of complexes)}\\&=-H_h(\delta_4(b))\quad\text{(follows from the definition of the map $H_h$)}\\&=\N_{\nu}(h)\delta_{3}(a)\quad\text{(because $-\delta_4(b)\in \mc L(\delta_3(a))$)}.
        \end{align*}
%\end{equation*}
\end{description}

 To prove that $\N_\nu(h)$ is a lifted reciprocity map it remains to show that it vanishes on elements of the form $a\wedge c, a\in \Lambda^2 \oL_\nu^\t, c\in k^\t$. Let  $b\in\mathcal L(a)$. Then $b\wedge c\in \mc L(a\wedge c)$ and we have:
 \begin{align*}
 & \N_\nu(h)(a\wedge c)=H_h(b\wedge c)\quad\text{(because $b\wedge c\in \mc L(a\wedge c)$)}  \\
 & =-\sum\limits_{\mu\in\dval(L)_{sp}}\N_{\mu}(h)\ts[4]_{\mu}(b\wedge c)\quad\text{(by the definition of $H_h$)}\\
 & =-\sum\limits_{\mu\in\dval(L)_{sp}}\N_{\mu}(h)(\ts[3]_{\mu}(b)\wedge c)\quad\text{(by the property of $\ts[4]_\mu$)}\\
 & =0\quad\text{(because $\N_\mu(h)$ is a lifted reciprocity map)}.
     \end{align*}
So we have proved that $\N_\nu(h)$ is a lifted reciprocity map.
\end{proof}

\subsection{Property of $\N_\nu$ under extensions of scalars}
\label{sec:norm_map:extensions_of_scalars}
The goal of this section is to prove the following statement:

\begin{proposition}
\label{prop:theta_extension_of_scalars}
Let $j_0\colon F\emb K$ be an embedding from $\F_1$ and $\nu\in \dval(F(t))$. Denote by $n$ the degree $[K:F]$ and by $j\colon F(t)\emb K(t)$ the unique extension of $j_0$ satisfying $j(t)=t$. For any $h\in \SRL(K)$, the following formula holds:
\begin{equation}
    \label{formula:extension_theta}
(\N_\nu\circ \SRL(j_0))(h)=\dfrac 1n\sum\limits_{\nu'\in \ext(\nu, K(t))}e_{\nn}f_{\nn}(\SRL(j_\nn)\circ \N_{\nu'})(h).
\end{equation}
\end{proposition}

This proposition is similar to Lemma 1.9 from \cite{suslin1979reciprocity}. (See also exercise {\MakeUppercase{\romannumeral 3}}.7.7 from \cite{KBook}).

\begin{definition}
\label{def:wt_theta}
Denote the right hand-side of formula (\ref{formula:extension_theta}) by $\wt \N_\nu(h)$. 
\end{definition}

Set $h_0=\SRL(j_0)(h)$. We need to show that $\N_\nu(h_0)=\wt \N_\nu(h)$. To do this, we need the following lemma:

\begin{lemma}
\label{Lemma:wt_theta_properties}
The following statements hold:
\begin{enumerate}
    \item For any $\nu\in \dval(F(t))$ the map  $\wt \N_\nu(h)$ is a lifted reciprocity map on the field $\ol{F(t)}_\nu$.
    \item Let $\wt h\in\SRL(F)$. For any $b\in \Lambda^4 F(t)^\times$ and all but finitely many $\nu\in \dval(F(t))$ we have $\N_\nu(\wt h)(\ts[4]_\nu(b))=0$ and moreover
$$\sum\limits_{\nu\in\dval(F(t))}\N_\nu(\wt h)(\ts[4]_\nu(b))=0.$$
\item For any $b\in \Lambda^4 F(t)^\times$ and all but finitely many $\nu\in \dval(F(t))$ we have $\wt\N_\nu(h)(\ts[4]_\nu(b))=0$ and moreover
$$\sum\limits_{\nu\in\dval(F(t))}\wt\N_\nu(h)(\ts[4]_\nu(b))=0.$$
\end{enumerate}
\end{lemma}

\begin{proof}[The deduction of Proposition \ref{prop:theta_extension_of_scalars} from Lemma \ref{Lemma:wt_theta_properties}]
We will prove the statement in two steps: first, we will check it when $\nu$ is special, then we reduce the case when $\nu$ is general to the previous case using Lemma \ref{Lemma:wt_theta_properties}.

\begin{enumerate}
    \item Let $\nu=\nu_{\infty, F}$. We recall that this is the valuation associated to the point $\infty\in\mathbb P^1_F$. In this case it is easy to see that the set $\ext(\nu,K(t))$ consists of only $1$ element, namely $\nu'=\nu_{\infty, K}$. We have $f_\nn=n, e_\nn=1$. Let us identify $\ol{F(t)}_{\nu_{\infty, F}}$ with $F$ and $\ol{K(t)}_{\nu_{\infty, K}}$ with $K$. Then the map $j_\nn$ is identified with $j_0$. Now, Definition \ref{def:wt_theta} gives: $\wt \N_\nu(h)=\SRL(j_\nn)(\N_{\nu'}(h))=\SRL(j_0)(h)=\N_\nu(\SRL(j_0)(h))=\mathcal N_\nu(h_0)$. Here we have used the definition of $\N_\nu$ and $\N_{\nu'}$ on special valuations (see the previous section). 
    \item Let $\nu$ be a special valuation different from $\nu_{\infty, F}$. By the item $(i)$ of the previous Lemma $\wt\N_\nu(h)$ is a lifted reciprocity map on the field $\ol{F(t)}_\nu$. Since $\ol{F(t)}_\nu\cong k(t)$, by Proposition \ref{prop:SRL_P1} any two lifted reciprocity maps on the field $\ol{F(t)}_\nu$ are equal. So $\N_\nu(h_0)=\wt\N_\nu(h)$.
    \item Let $\nu$ be a general valuation. We need to show that for any $a\in\Lambda^3\ol{F(t)}_\nu^\times$, the following formula holds: $\N_\nu(h_0)(a)=\wt\N_\nu(h)(a)$. Choose some $b\in \mathcal L(a)$. We have: 
    \begin{align*}
    &\N_\nu(h_0)(a)=\\
& =\sum\limits_{\nu\in\dval(F(t)_{gen})}\N_\nu(h_0)(\ts[4]_\nu(b))\quad\text{(Because $b\in\mathcal L(a)$)}\\
&=-\sum\limits_{\nu\in\dval(F(t)_{sp})}\N_\nu(h_0)(\ts[4]_\nu(b))\quad\text{(By item (ii) of Lemma \ref{Lemma:wt_theta_properties})}\\
&=-\sum\limits_{\nu\in\dval(F(t)_{sp})}\wt\N_\nu(h)(\ts[4]_\nu(b))\quad\text{(By item (i) and (ii) of this proof)}\\
&=\sum\limits_{\nu\in\dval(F(t)_{gen})}\wt\N_\nu(h)(\ts[4]_\nu(b))\quad\text{(By item (iii) of Lemma \ref{Lemma:wt_theta_properties})}\\
&=\wt\N_\nu(h)(a)\quad\text{(Because $b\in\mathcal L(a)$).}
     \end{align*}
    So $\N_\nu(h_0)=\wt\N_\nu(h)$.
\end{enumerate}
\end{proof}

\begin{proof}[The proof of Lemma \ref{Lemma:wt_theta_properties}]
\begin{enumerate}
    \item The set of all lifted reciprocity maps on some field has a structure of an affine set over $\mathbb Q$, see Remark \ref{rem:SRL_affine}. This means that if $X$ is a finite set and $h_{\mu}, \mu\in X$ are some lifted reciprocity maps on the field $\ol{F(t)}_\nu$, then for any $\alpha_\mu\in \mathbb Q$ satisfying  $\sum\limits_{\mu\in X}\alpha_\mu=1$ the map defined by the formula  $\sum\limits_{\mu\in X}\alpha_\mu h_\mu$ is a lifted reciprocity map on the field $\ol{F(t)}_\nu$. Applying this statement to $X=\ext(\nu,K(t)), \alpha_\mu=\dfrac{e_{\mu|\nu}f_{\mu|\nu}}{[K:F]}$ and $h_\mu=\SRL(j_{\mu|\nu})(\N_\mu(h))$, we get the statement of the lemma. (The formula $\sum\limits_{\mu\in X}\alpha_\mu=1$ follows from formula (\ref{formula:degree_as_sum})).
    \item The first statement follows from Lemma \ref{lemma:finitnes_of_sum}. Let us prove the second statement. It follows from Theorem \ref{th:main_exact_sequence} that the sets $\mathcal L(a), a\in \Lambda^3 \oL_\mu^\times, \mu\in \dval(F(t))$ generate $\Lambda^4 F(t)^\times$ as a vector space. So we can assume that $b\in\mathcal L(a)$ for some $a\in \Lambda^3 \oL_\mu^\times, \mu\in \dval(F(t))$. We need to show the following equality:
$$\sum\limits_{\nu\in\dval(F(t))_{gen}}\N_\nu(\wt h)(\ts[4]_\nu(b))=-\sum\limits_{\nu\in\dval(F(t))_{sp}}\N_\nu(\wt h)(\ts[4]_\nu(b)).$$
By the definition of the set $\mathcal L(a)$ (see Definition \ref{def:lift}) the left hand-side is equal to $\N_{\mu}(\wt h)(a)$. On the other hand the right hand side is equal to $H_{\wt h}(b)$ which is exactly the definition of $\N_{\mu}(\wt h)(a)$ (see Definition \ref{def:theta_gen}). 
\item The first statement of the lemma follows from Proposition \ref{lemma:finitnes_of_sum} and item $(i)$ of this lemma. (We recall that any lifted reciprocity map on the field $\ol{F(t)}_\nu$ is zero on the image of the multiplication map $\ol{F(t)}_\nu^\t\otimes \Lambda^2 k^\t\to \Lambda^3 \ol{F(t)}_\nu^\t$).

To prove the second statement, let us rewrite the element $\wt\N_\nu(h)(\ts[4]_\nu(b))$ as follows:
\begin{align*}
& \wt\N_\nu(h)(\ts[4]_\nu(b))\\
&=\dfrac 1n\sum\limits_{\nu'\in\ext(\nu, K(t))}e_{\nn} f_{\nn}\SRL(j_\nn)(\N_{\nu'}(h)(\ts[4]_\nu(b)))\quad\text{(By the definition of }\wt\N_\nu(h))\\
&=\dfrac 1n\sum\limits_{\nu'\in\ext(\nu, K(t))}e_{\nn}\N_{\nu'}(h)((j_\nn)_*(\ts[4]_\nu(b)))\quad\text{(By the definition of }\SRL(j_\nn))\\
&=\dfrac 1n \sum\limits_{\nu'\in\ext(\nu, K(t))}\N_{\nu'}(h)(\ts[4]_{\nu'}(j_*(b)))\quad\text{(By formula  \eqref{formula:functoriality_of_tame_symbol})}.
     \end{align*}
     
     So we have:
  \begin{align*}
& \sum\limits_{\nu\in\dval(F(t))}\wt\N_\nu(h)(\ts[4]_\nu(b))=\sum\limits_{\nu\in\dval(F(t))}\dfrac 1n\sum\limits_{\nu'\in\ext(\nu, K(t))}\N_{\nu'}(h)(\ts[4]_{\nu'}(j_*(b)))\\
&=\dfrac 1n\sum\limits_{\nu'\in\dval(K(t))}\N_{\nu'}(h)(\ts[4]_{\nu'}(j_*(b))).
     \end{align*}
The last expression is zero by item $(ii)$ of this lemma applied to the field $K$, the lifted reciprocity map $h\in\SRL(K)$  and the element $j_*(b)\in \Lambda^4 K(t).$
\end{enumerate}
\end{proof}

\subsection{The proof of Theorem \ref{th:norm_map}}
\label{sec:norm_map:finall_proof}
In this section we will use results of Section \ref{sec:norm_map:theta} to construct the norm map on lifted reciprocity maps. We follow ideas from \cite[\S 1]{suslin1979reciprocity} (see also \cite{MILNOR1969/70, bass1973milnor, kato1980generalization}). 

\begin{definition}
\label{def:norm_map}
Let $j\colon F\emb K$ be an extension of some fields from $\F_1$. Let $a$ be some generator of $K$ over $F$. Denote by $p_a\in F[t]$ the minimal polynomial of $a$ over $F$. Denote by $\nu_a$ the corresponding valuation. The residue field $\overline{F(t)}_{\nu_{a}}$ is canonically isomorphic to $K$. So we get the map $\N_{{\nu_a}}\colon\SRL(F)\to \SRL(K)$, which we denote by $N_{K/F, a}$. This map is called \emph{the norm map}. 
\end{definition}

\begin{proof}[The proof of Theorem \ref{th:norm_map}]
The proof of Theorem \ref{th:norm_map} goes as follows. First of all we will show that the map $N_{K/F, a}$ is well behaved with respect to extension of scalars (Lemma \ref{lemma:functoriality_of_norms}). This will follow directly from Proposition \ref{prop:theta_extension_of_scalars}. Then we will prove Lemma \ref{lemma:SRL_is_surjective} stating that $N_{K/F, a}$ is a right inverse for $\SRL(j)$. This will show that $\SRL(j)$ is surjective and that in the case $K=F$ the map $N_{K/F,a}$ does not depend on $a$. Then we will prove Proposition \ref{prop:norm_does_not_depend_on_a} stating that the map $N_{K/F, a}$ does not depend on $a$. By that moment, for any field extension $F\emb K$ we will have constructed the canonical norm map $N_{K/F}$ and will have proved item $(i)$ of Theorem \ref{th:norm_map}. The item $(ii)$ will follow from Proposition \ref{prop:tower}.
\end{proof}

\begin{lemma}
\label{lemma:functoriality_of_norms}
Let $j\colon F_1\emb K,  F_1\emb F_2, $ be extensions and $F_2\otimes_{F_1} K=\bigoplus\limits_{i=1}^m F_{2,i}$. Denote by $j_i$ the natural embedding $F_2\emb F_{2,i}$. Let $n=[K:F_1]$ and $n_i=[F_{2,i}:F_2]$. Let $a$ be a generator of $F_2$ over $F_1$. Denote by $a_i$ the corresponding generators of $F_{2,i}$ over $K$. The following diagram is commutative:
     \begin{equation}
    \label{diagram:functoriality_of_norm}
        \begin{tikzcd}[row sep=huge,column sep=huge]
	{\SRL(F_1)} & {\SRL(F_2)} \\
	{\SRL(K)} & {\bigoplus\limits_{i=1}^m\SRL(F_{2,i}).}
	\arrow["N_{F_2/F_1,a}",from=1-1, to=1-2]
	\arrow["(N_{F_{2,i}/K,a_i})",from=2-1, to=2-2]
	\arrow["\SRL(j)",from=2-1, to=1-1]
	\arrow["\sum\limits_{i=1}^m \dfrac{n_i}{n}\SRL(j_i)"',from=2-2, to=1-2]
\end{tikzcd}
    \end{equation}
\end{lemma}

\begin{proof}
It follows from Proposition \ref{prop:theta_extension_of_scalars} that for any $\nu\in \dval(F_1(t))$ the following diagram is commutative:
 \begin{equation}
 \label{diagram:functoriality_of_norm_valuations}
        \begin{tikzcd}[row sep=huge,column sep=huge]
	{\SRL(F_1)} & {\SRL(\overline{F_1(t)}_\nu)} \\
	{\SRL(K)} & {\bigoplus\limits_{\nu'\in \ext(\nu, K(t))}\SRL(\overline{K(t)_{\nu'}}).}
	\arrow["\N_{\nu}",from=1-1, to=1-2]
	\arrow["(\N_{\nu'})",from=2-1, to=2-2]
	\arrow["\SRL(j)",from=2-1, to=1-1]
	\arrow["\sum\limits_{\nu'\in \ext(\nu,K(t))} \dfrac{e_{\nn}f_{\nn}}n\SRL(j_{\nn})"',from=2-2, to=1-2]
\end{tikzcd}
    \end{equation}

Let us apply this statement in the case when $\nu$ is equal to $\nu_a$. Let $p_a=\prod\limits_{i}^mp_{a,i}$ be the decomposition of $p_a$ in the field $K(t)$. The set $\ext(\nu, K(t))$ is in bijection with the irreducible factors of $p_a$ in $K(t)$. Denote by $\nu_i\in \ext(\nu, K(t))$ the valuation corresponding  to $p_{a,i}$. We have $F_{2,i}\cong \overline{K(t)}_{\nu_i}$. The embeddings $j_i$ correspond to the embeddings $j_{\nu_i|\nu}$. Since the polynomial $p_a$ is separable, we have $f_{\nn}=1$, and so $e_{\nn}f_{\nn}=[F_{2,i}:F_2]$. So the diagram (\ref{diagram:functoriality_of_norm_valuations}) can be identified with (\ref{diagram:functoriality_of_norm}).

\end{proof}

\begin{lemma}
\label{lemma:SRL_is_surjective}
For any embedding $j\colon F_1\emb F_2$, we have $\SRL(j)\circ N_{F_2/F_1, a}=id$. In particular, the map $\SRL(j)$ is surjective and in the case $F_1=F_2$, the map $N_{F_2/F_1, a}$ is the identity map.
\end{lemma}

\begin{proof}Let $n=[F_2:F_1]$. We need to show that for any $h\in\SRL(F_1)$ and $x\in\Lambda^3 F_1^\times$ the following formula holds: $N_{F_2/F_1, a}(h)(j_*(x))=n\cdot h(x)$.
 Consider the element $b=p_a\wedge x\in \Lambda^{4} F_1(t)^\t$, where $p_a$ is the minimal polynomial of $a$ over $F_1$.  By item $(ii)$ of Lemma \ref{Lemma:wt_theta_properties} we have:
    $$\sum\limits_{\nu\in \dval(F_1(t))_{gen}}\N_\nu(h)(\ts[4]_\nu(b))+\sum\limits_{\nu\in \dval(F_1(t))_{sp}}\N_\nu(h)(\ts[4]_\nu(b))=0.$$
    We have $\ts[4]_{\nu_{p_a}}(b)=x$ and this is the only general valuation satisfying $\ts[4]_\nu(b)\ne 0$. So the first term is equal to $\N_{\nu_{p_a}}(h)(x)$ which is equal to $N_{F_2/F_1, a}(h)(x)$. So we have:
    $$N_{F_2/F_1,a}(h)(x)=-\sum\limits_{\nu\in \dval(F_1(t))_{sp}}\N_\nu(h)\ts[4]_\nu(b).$$
    It is easy to see that there is only one special valuation $\nu$ such that $\N_\nu(h)(\ts[4]_\nu(b))\ne 0$, namely $\nu_{F_1,\infty}$.  We have  $\ts[4]_{\nu_{F_1,\infty}}(b)=-nx$. Since $\N_{\nu_{F_1,\infty}}(h)$ can be identified with $h$, we have $N_{F_2/F_1,a}(h)(x)=-\N_{\nu_{F_1,\infty}}(h)(\ts[4]_{\nu_{F_1,\infty}}(b))=n\cdot h(x)$.
\end{proof}

\begin{proposition}
\label{prop:norm_does_not_depend_on_a}
The map $N_{F_2/F_1,a}$ does not depend on $a$. 
\end{proposition}

We denote the map $N_{F_2/F_1,a}$ simply by $N_{F_2/F_1}$.

\begin{proof}
Let $j\colon F_1\emb K$ be a field extension of $F_1$ satisfying $F_2\otimes_{F_1} K\cong K^{\oplus[F_2:F_1]}$.  We apply Lemma \ref{lemma:functoriality_of_norms}. By definition of $K$ for any $n$ we have $F_{2,i}\cong K$. By Lemma \ref{lemma:SRL_is_surjective} the maps $N_{F_{2,i}/K, a_i}$ are the identity maps. We conclude that in the diagram from Lemma \ref{lemma:functoriality_of_norms} all the maps except maybe $N_{F_2/F_1,a}$ do not depend on $a$. So the map $N_{F_2/F_1, a}$ does not depend on $a$ on the image of $\SRL(j)$. By the previous lemma this image coincides with $\SRL(F_1)$.
\end{proof}

\begin{proposition}
\label{prop:tower}
If $F_1\subset F_2\subset F_3$ is a tower of extensions from $\F_1$ then $N_{F_3/F_1}=N_{F_3/F_2}\circ N_{F_2/F_1}$.
\end{proposition}

\begin{proof}
Let $j\colon F_1\emb K$ be a field extension. Denote $F_2\otimes_{F_1}K\cong \bigoplus_{i=1}^{n_{2}} F_{2,i}$ and $F_{3}\otimes_{F_2} F_{2,i} \cong \bigoplus_{s=1}^{n_{3,i}}F_{3,i,s}$. By associativity of tensor product we have $F_3\otimes_{F_1} K\cong \bigoplus\limits_{i,s}F_{3,i,s}$ Denote by $j_{i,s}$ the natural embeddings $F_3\emb F_{3,i,s}$. Let $n_{i,s}=[F_{3,i,s}:F_3]$. Let $n$ be the degree of $F_3$ over $F_1$.
Repeated application of Lemma \ref{lemma:functoriality_of_norms} together with Proposition \ref{prop:norm_does_not_depend_on_a} shows that  the following diagram is commutative:
 \begin{equation*}
        \begin{tikzcd}[row sep=huge,column sep=huge]
	{\SRL(F_1)} & & {\SRL(F_3)} \\
	{\SRL(K)} & & {\bigoplus\limits_{i,s}\SRL(F_{3,i,s}).}
	\arrow["N_{F_3/F_2}\circ N_{F_2/F_1}",from=1-1, to=1-3]
	\arrow["(N_{F_{3,i,s}/F_{2,i}}\circ N_{F_{2,i}/K})",from=2-1, to=2-3]
	\arrow["\SRL(j)",from=2-1, to=1-1]
	\arrow["\sum\limits_{i,s} \dfrac{n_{i,s}}{n}\SRL(j_{i,s})"',from=2-3, to=1-3]
\end{tikzcd}
    \end{equation*}

Choose $K$ such that $F_3\otimes_{F_1}K\cong K^{\oplus [F_3:F_1]}$. It follows that $F_{3,i,s}\cong F_{2,i}\cong K$. So the bottom maps in the above diagram are the identity maps. Let us compare this diagram with diagram (\ref{diagram:functoriality_of_norm}) for $F_2=F_3$. We see that the left, right and bottom maps are the same. Since $\SRL(j)$ is surjective, the statement of the proposition follows.
\end{proof}

%Theorem \ref{th:norm_map} follows directly from Lemma \ref{lemma:SRL_is_surjective},  Proposition \ref{prop:norm_does_not_depend_on_a}, Proposition \ref{prop:tower} and Proposition \ref{prop:Ha=Hb}.
\section{The proofs of the main results}
\label{sec:finall_proofs}

\subsection{The proof of Theorem \ref{th:limit_of_functor}}
\label{sec:proof_of_main_result}
Let $F_0$ be the field $k(x)$ and $L=F_0(y)$. For any $\nu\in\dval(L)$ denote by $\sigma_\nu\in \SRL(\oL_\nu)$ the lifted reciprocity map given by the formula $\N_\nu(\mathcal H_{k(x)})$. ($\mathcal H_{k(x)}$ is a lifted reciprocity map from Proposition \ref{prop:SRL_P1}). Denote by $\lambda$ the involution of $L$ fixing $k$ and interchanging $x$ and $y$. The map $\lambda$ induces the natural map $\dval(L)\to\dval(L)$ given by the formula $\lambda(\nu)(f)=\nu(\lambda(f))$. Denote by $\ol \lambda_{\nu}$ the natural map $\ol\lambda_\nu\colon \oL_\nu\to \oL_{\lambda(\nu)}$.

\begin{lemma}
For any $\nu\in \dval(L)$ we have $\sigma_\nu=\SRL(\ol\lambda_\nu)(\sigma_{\lambda(\nu)})$.
\end{lemma}

\begin{proof}
For any $\nu\in\dval(L)$ denote $\wt\sigma_\nu=\SRL(\ol\lambda_\nu)(\sigma_{\lambda(\nu)})$. We need to show that $\sigma_\nu=\wt\sigma_\nu$. 
\begin{enumerate}
    \item Let us assume that $\nu$ is special. In this case $\oL_\nu\cong k(t)$ and the statement follows from Proposition \ref{prop:SRL_P1}.
    \item Let $\nu$ be general. It is easy to see that $\wt \sigma$ satisfies item $(ii)$ of Lemma \ref{Lemma:wt_theta_properties}. Namely, for any $b\in\Lambda^4 L^\t$ and for all but finitely many $\nu\in\dval(L)$ we have $\wt\sigma_\nu(\ts[4]_\nu(b))=0$ and moreover the following formula holds:
    $$\sum\limits_{\nu\in\dval(L)}\wt\sigma_\nu(\ts[4]_\nu(b))=0.$$
    Now, the case of a general valuation can be reduced to the case of a special one in the same way as it was done in the proof of Proposition \ref{prop:theta_extension_of_scalars}.
\end{enumerate}
\end{proof}

\begin{proposition}
\label{prop:Ha=Hb}
For any $a,b\in F\bs k$ we have $$N_{F/k(a)}(\mc H_{k(a)})=N_{F/k(b)}(\mc H_{k(b)}).$$ In particular, the element $\mc H_F:=N_{F/k(a)}(\mc H_{k(a)})$ does not depend on the choice of $a\in F\bs k$.
\end{proposition}

\begin{proof} We first prove the statement in the case when $a$ and $b$ generate $F$ over $k$.

Let us realize the field $F$ as filed of fractions on some smooth projective curve $X$. The elements $a,b\in F$ induce the maps $\varphi_a, \varphi_b\colon X\to \mathbb P^1$. Define the maps $\psi_{a,b}, \psi_{b,a}\colon X\to \mathbb P^1\t \mathbb P^1$ given by the formulas $\psi_{a,b}=(\varphi_a, \varphi_b), \psi_{b,a}=(\varphi_b, \varphi_a)$.  Let $S=\mathbb P^1\times \mathbb P^1$. We can identify $k(S)$ with $L=k(x)(y)$. Denote by $\mu_{a,b}$ (resp. $\mu_{b,a}$) the valuation of $L$ corresponding to the image of $\psi_{a,b}$ (resp. $\psi_{b,a}$). We have $\lambda(\mu_{a,b})=\mu_{b,a}$. Denote by $\theta_{a,b}$ (resp. $\theta_{b,a}$) the canonical isomorphism $F\to \overline{L}_{\nu_{a,b}}$ (resp. $F\to \overline{L}_{\nu_{b,a}}$).
We recall that the definition of $\sigma_\nu$ was given in the beginning of this section.
By the definition of the norm map, we have:
$$N_{F/k(a)}(\mc H_{k(a)})=\SRL(\theta_{a,b})(\sigma_{\nu_{a,b}}), N_{F/k(b)}(\mc H_{k(b)})=\SRL(\theta_{b,a})(\sigma_{\nu_{b,a}}).$$
So we need to show that $\SRL(\theta_{a,b})(\sigma_{\nu_{a,b}})=\SRL(\theta_{b,a})(\sigma_{\nu_{b,a}}).$
Denote the map $\ol \lambda_{\mu_{a,b}}\colon \oL_{\mu_{a,b}}\to\oL_{\mu_{b,a}}$ simply by $\ol \lambda$. We have $\ol\lambda\circ\theta_{a,b}=\theta_{b,a}$.   
 We get:
\begin{equation*}
     \begin{split}
         \SRL(\theta_{b,a})(\sigma_{\mu_{b,a}})=\SRL(\overline {\lambda}\circ\theta_{a,b})(\sigma_{\mu_{b,a}})=
         (\SRL(\theta_{a,b})\circ\SRL(\overline{\lambda}))(\sigma_{\mu_{b,a}})=\\=\SRL(\theta_{a,b})(\SRL(\overline{\lambda})(\sigma_{\mu_{b,a}}))=\SRL(\theta_{a,b})(\sigma_{\mu_{a,b}}).
     \end{split}
 \end{equation*}
 In the last formula we have used the previous lemma.
 
So we have proved that for any $a,b\in F\bs k$ generating $F$ over $k$ we have $$N_{F/k(a)}(\mc H_{k(a)})=N_{F/k(b)}(\mc H_{k(b)}).$$ Now the first statement of the proposition follows from the following fact: for any $a,b\in F\bs k$, there is $c\in F\bs k$ such that the pairs $(a,c),(b,c)$ generate $F$ over $k$. 

The second statement is a reformulation of the first statement.

\end{proof}
We recall that we work $\mathbb Q$-linearly. The following lemma is well-known:
\begin{lemma}
\label{lemma:Galois_descent}
Let $j\colon F_1\emb F_2$ be a Galois extension with the Galois group equal to $G$. The natural map $j_*\colon K_n^M(F_1)\to K_n^M(F_2)$ induces an isomorphism $K_n^M(F_1)\xrightarrow{\sim} K_n^M(F_2)^G$
\end{lemma}

\begin{proof}
There is the norm map $\wh N_{F_2/F_1}\colon K_n^M(F_2)\to K_n^M(F_1)$,  satisfying the following properties (see \cite{KBook}):
\begin{enumerate}
    \item The composition $\wh N_{F_2/F_1}\circ j_*\colon K_n^M(F_1)\to K_n^M(F_1)$ is the multiplication by the integer $[F_2:F_1]$.
    \item The composition $j_*\circ \wh N_{F_2/F_1}$ is equal to $\sum\limits_{g\in G}g_*$.
    These statements imply that the kernel and the cokernel of the natural map $K_n^M(F_1)\to K_n^M(F_2)^G$ is annihilated by the multiplication on $[F_2:F_1]$. Therefore the map $K_n^M(F_1)\to K_n^M(F_2)^G$ is a rational isomorphism.
\end{enumerate}
\end{proof}

\label{sec:proof_of_main_resulats}

\begin{proof}[The proof of Theorem \ref{th:limit_of_functor}]

\begin{description}
\item[Existence] Let $F\in \F_1$. Choose some embedding $$j\colon k(t)\emb F.$$ Define the element $\mc H_F$ by the formula $\mc H_F:=N_{F/k(t)}(\mc H_{k(t)}).$ By Proposition \ref{prop:Ha=Hb} this element does not depend on $j$. 
We need to show that if $j'\colon F_1\to  F_2$ is an embedding, then $$\SRL(j')(\mathcal H_{F_2})=\mc H_{F_1}.$$ This follows from the properties of the norm map (see Theorem \ref{th:norm_map}):
\begin{equation*}
    \begin{split}
        \SRL(j')(\mc H_{F_2})=\SRL(j')N_{F_2/k(t)}(\mc H_{k(t)})=\SRL(j')N_{F_2/F_1}N_{F_1/k(t)}(\mc H_{k(t)})=\\
        =(\SRL(j')\circ N_{F_2/F_1})(N_{F_1/k(t)}\mc H_{k(t)})=\mc H_{F_1}.
    \end{split}
\end{equation*}
\item[Uniqueness]
Let $\mc H_F, \mc H_F', F\in\F_1$ be two families of lifted reciprocity maps such that for any $j\colon F_1\emb F_2$ we have $\SRL(j)(\mc H_{F_2})=\mc H_{F_1}$ and $\SRL(j)(\mc H'_{F_2})=\mc H'_{F_1}$. We need to show that $\mc H_F=\mc H_F'$ for any $F\in \F_1$. By Proposition \ref{prop:SRL_P1}, this is true when $F=k(t)$. Let $F$ be any field. There is a field $F'\in\F_1$ together with two embeddings $j_1\colon F\emb F', k(t) \emb F'$ such that $F'/k(t)$ is Galois. Since $\mc H_F=\SRL(j_1)(\mc H_{F'})$ and $\mc H_F'=\SRL(j_1)(\mc H_{F'}')$, it is enough to prove the statement for $F'$. Denote by $G$ the Galois group of $F'$ over $k(t)$. Since $\mc H_{F'}$ and $\mc H_{F'}'$ are invariant under the group $G$, it is enough to prove that they are equal on the subgroup $\left(\Lambda^{3} F'^\t\right)^G$. By Lemma \ref{lemma:Galois_descent}, we have $\left(K_{3}^M(F')\right)^G=K_{3}^M(k(t))$. It follows that $\left(\Lambda^{3} F'^\t\right)^G$ is generated by the image of $\delta_3$ and by the elements coming from $k(t)$. On the image of $\delta_3$ the maps $\mc H_{F'}$ and $\mc H_{F'}'$ coincide because they are lifted reciprocity maps. On the elements coming from $k(t)$ they coincide because $\mc H_{k(t)}=\mc H_{k(t)}'$.

\end{description}
\end{proof}

\subsection{The proof of Corollary \ref{cor:reciprocity_map_for_surfaces}}
\label{sec:proof_of_reciprocity_law_for_surfaces}

Let $L\in \F_2$. Define the map $H_L\colon \Lambda^{4} L^\t\to {\P}(k)$ by the formula:

$$H_L(b)=\sum\limits_{\nu\in \dval(L)}\mc H_{\oL_\nu}(\ts[4]_\nu(b)).$$

This formula is well-defined by Lemma \ref{lemma:finitnes_of_sum}.  The following lemma is a consequence of Theorem \ref{th:limit_of_functor}:
\begin{lemma}
\label{lemma:functoriality_for_H}
If $j\colon L\hookrightarrow M$ is an extension of some fields from $\F_2$ then for any $b\in \Lambda^{4} L^\t$ we have $H_{L}(b)=\dfrac 1{[M:L]}H_{M}(j_*(b))$.
\end{lemma}

\begin{proof}
Let $\nu\in \dval(L)$. It is enough to show  the following formula:
$$\mc H_{\oL_\nu}\ts[4]_\nu(b)=\dfrac 1{[M:L]}\sum\limits_{\nu'\in \ext(\nu,M)}\mc H_{\oM_{\nu'}}\ts[4]_{\nu'}(j_*(b)). $$
By formula (\ref{formula:functoriality_of_tame_symbol}), we have: $$\ts[4]_{\nu'}(j_*(b))=e_{\nn} \cdot j_{\nn}(\ts[4]_\nu(b)).$$
If $s\colon F_1\emb F_2$ is an arbitrary extension from $\F_1$ then for any $x\in\Lambda^3 F_1^\times$, we have: $$\mathcal H_{F_2}({s_*(x)})=[F_2:F_1]\cdot \mathcal H_{F_1}(x).$$ Applying this in the case $s=j_{\nn}\colon \overline L_\nu\emb \overline M_{\nu'}$, we get: $$\mc H_{\oM_{\nu'}}((j_{\nn})_*(\ts[4]_\nu(b)))=f_{\nn} \mc H_{\oL_\nu}(\ts[4]_\nu(b)).$$ 
So we get:
$$\mc H_{\oM_{\nu'}}(\ts[4]_{\nu'}(j_*(b)))=e_{\nn}f_{\nn} \cdot\mc H_{\oL_\nu}(\ts[4]_\nu(b)).$$
Using this formula, we obtain:
\begin{align*}
   &\dfrac 1{[M:L]}\sum\limits_{\nu'\in \ext(\nu, M)} \mc H_{\oM_{\nu'}}(\ts[4]_{\nu'}j_*(b))=\\
   &=\dfrac 1{[M:L]}\sum\limits_{\nu'\in \ext(\nu, M)}e_{\nn}f_{\nn} \mc H_{\oL_\nu}(\ts[4]_\nu(b))= \\&=\mc H_{\oL_\nu}(\ts[4]_\nu(b))\dfrac 1{[M:L]}\sum\limits_{\nu'\in \ext(\nu,M)}e_{\nn}f_{\nn}=H_{\oL_\nu}(\ts[4]_\nu(b)).
   \end{align*}
The last equality follows from the formula $\sum\limits_{\nu'\in \ext(\nu, M)}e_{\nn}f_{\nn}=[M:L]$.
\end{proof}

\begin{proof}[Proof of Corollary \ref{cor:reciprocity_map_for_surfaces}]
We need to show that $H_L=0$ for any $L\in \F_2$. Let us prove that it is true when $L=k(x)(y)$. 

Let us represent the field $L$ as $F_0(y)$, where $F_0=k(x)$. For any $\nu\in \dval(F_0(y))$ define $\sigma_\nu=\mathcal N_\nu(\mathcal H_{k(x)})\in \SRL(\oL_\nu)$. Let us prove that $\sigma_\nu=\mathcal H_{\oL_\nu}$. When $\nu$ is special it is true by Proposition \ref{prop:SRL_P1}. When $\nu$ is general it follows from the definition of $\mc H_{\overline L_\nu}$. Now the formula
$$\sum\limits_{\nu\in\dval(L)}\mathcal H_{\oL_\nu}(\ts[4]_\nu(b))$$
follows from item $(ii)$ of Lemma \ref{Lemma:wt_theta_properties}.

Let us prove the statement for an arbitrary $L$. There is a field $L'$ together with two finite extensions $$j\colon k(x)(y)\emb L', j'\colon L\emb L'$$ such that $j$ is a Galois extension. Lemma \ref{lemma:functoriality_for_H} shows that it is enough to prove the statement for $L'$. Denote the Galois group of $j$ by $G$. Since $H_{L'}$ is invariant under $G$, it is enough to prove that $H_{L'}$ is zero on the subgroup $\left(\Lambda^4 L'^\t\right)^G$. By Lemma \ref{lemma:Galois_descent}, we have $\left(K_{4}^M(L')\right)^G=K_{4}^M(k(x)(y))$ and so the group $\left(\Lambda^{4} L'^\t\right)^G$ is generated by the image of $\delta_4$ and by the elements coming from $k(x)(y)$. Vanishing of $H_{L'}$ on the elements coming from $k(x)(y)$ follows from Lemma \ref{lemma:functoriality_for_H} together with the formula $H_{k(x)(y)}=0$. Let us prove that  $H_{L'}$ is zero on the image of the map $\delta_4\colon \Gamma(L',4)_3\to \Gamma(L',4)_4$. For any $b\in \Gamma(L',4)_3$ we have
\begin{equation*}
\begin{split}
    H_{L'}(\delta_4(b))=\sum\limits_{\nu\in \dval(L')}\mc H_{\overline {L'}_\nu}\ts[4]_\nu\delta_4(b)=-\sum\limits_{\nu\in \dval(L')}\mc H_{\overline {L'}_\nu}\delta_3\ts[4]_\nu(b)=\\=-\sum\limits_{\nu\in \dval(L')}\sum\limits_{\nu'\in \val(\overline {L'}_\nu)}\ts_{\nu'}\ts[4]_{\nu}(b)=0.
    \end{split}
\end{equation*}
Here the second equality is true because $\ts[4]_\nu$ is morphism of complexes, the third equality is true because $\mc H_{\overline {L'}_\nu}$ is a lifted reciprocity map and the fourth equality follows from Theorem \ref{th:Parshin_sum_point_curve}.
\end{proof}

\bibliographystyle{alpha}  
\bibliography{mylib} 

\end{document}